\documentclass[
	english
]{scrartcl}

\usepackage[T1]{fontenc}
\usepackage[utf8]{inputenc}

\usepackage{amsmath,amsfonts,amsthm,amssymb}
\usepackage[colorlinks,citecolor=blue]{hyperref}
\usepackage{doi}
\usepackage{cancel}
\usepackage{relsize}
\usepackage{booktabs}

\usepackage{changes}

\usepackage[shortlabels]{enumitem}
\setlist[enumerate]{label=(\alph*)}

\usepackage{tikz,pgfplots,pgfplotstable}
\usetikzlibrary{arrows.meta}
\usetikzlibrary{decorations.markings}

\usepackage[figure]{hypcap}
\usepackage{graphicx}
\graphicspath{{./img/}}
\usepackage{subcaption} 

\usepackage{preamble}

\usepackage{marginnote}

\usepackage[section]{algorithm}
\usepackage{algpseudocode}
\usepackage{varwidth}

\makeatletter
\newcommand{\algmargin}{\the\ALG@thistlm}
\makeatother
\algnewcommand{\parState}[1]{\State\parbox[t]{\dimexpr\linewidth-\algmargin}{\strut #1\strut}}
\usepackage[capitalise,nameinlink]{cleveref}

\crefname{figure}{Figure}{Figures}
\crefname{table}{Table}{Tables}

\allowdisplaybreaks

\numberwithin{equation}{section}
\counterwithin{figure}{section}

\newcommand\norm[1]{\left\Vert#1\right\Vert}
\newcommand\nnorm[1]{\Vert#1\Vert}

\newcommand\N{\mathbb{N}}
\newcommand\R{\mathbb{R}}

\newcommand\ff{\mathsf{f}}
\renewcommand\gg{\mathsf{g}}
\newcommand\hh{\mathsf{h}}

\newcommand\tto{\rightrightarrows}

\newcommand{\dist}{\operatorname{dist}}

\newcommand{\dom}{\operatorname{dom}}

\newcommand{\gph}{\operatorname{gph}}
\newcommand{\epi}{\operatorname{epi}}

\DeclareMathOperator*{\argmin}{\operatorname{argmin}}

\DeclareMathAlphabet{\mathpzc}{OT1}{pzc}{m}{it}
\newcommand\oo{\mathpzc{o}}

\DeclareFontFamily{U}{mathx}{\hyphenchar\font45}
\DeclareFontShape{U}{mathx}{m}{n}{
      <5> <6> <7> <8> <9> <10>
      <10.95> <12> <14.4> <17.28> <20.74> <24.88>
      mathx10
      }{}
\DeclareSymbolFont{mathx}{U}{mathx}{m}{n}
\DeclareFontSubstitution{U}{mathx}{m}{n}
\DeclareMathAccent{\widecheck}{0}{mathx}{"71}
\DeclareMathAccent{\wideparen}{0}{mathx}{"75}

\newtheorem{theorem}{Theorem}[section]
\newtheorem{lemma}[theorem]{Lemma}
\newtheorem{proposition}[theorem]{Proposition}
\newtheorem{assumption}[theorem]{Assumption}

\newtheorem{corollary}[theorem]{Corollary}
\newtheorem{remark}[theorem]{Remark}
\newtheorem{definition}[theorem]{Definition}
\newtheorem{example}[theorem]{Example}

\counterwithin{table}{section}

\usetikzlibrary{arrows}

\tikzset{
    punkt/.style={
           rectangle,
           draw=white, very thick,
           text width=9em,
           minimum height=1.5em,
           text centered}
}

\begin{document}

\title{Revisiting implicit variables in mathematical optimization: simplified modeling and a numerical evidence}
\author{%
	Patrick Mehlitz%
	\footnote{%
		Philipps-Universität Marburg,
		Department of Mathematics and Computer Science,
		35032 Marburg,
		Germany,
		\email{mehlitz@uni-marburg.de},
		\orcid{0000-0002-9355-850X}%
		}
}

\publishers{}
\maketitle

\begin{abstract}
	Implicit variables of an optimization problem are used 
	to model variationally challenging feasibility conditions in a tractable way
	while not entering the objective function.
	Hence, it is a standard approach to treat implicit variables as explicit ones.
	Recently, it has been shown in terms of a comparatively complex model problem
	that this approach, generally, is theoretically disadvantageous as the surrogate problem typically suffers
	from the presence of artificial stationary points and the need for stronger
	constraint qualifications.
	The purpose of the present paper is twofold.
	First, it introduces a much simpler and easier accessible model problem 
	which can be used to recapitulate and even broaden the aforementioned findings.
	Indeed, we will extend the analysis to two more classes of stationary points
	and the associated constraint qualifications.
	These theoretical results are accompanied by illustrative examples from
	cardinality-constrained, vanishing-constrained, and bilevel optimization.
	Second, the present paper illustrates, 
	in terms of cardinality-constrained portfolio optimization problems,
	that treating implicit variables as explicit ones may also be
	disadvantageous from a numerical point of view.
\end{abstract}


\begin{msc}	
	\mscLink{49J53}, \mscLink{90C30}, \mscLink{90C33}
\end{msc}

\section{Introduction}\label{sec:introduction}

A variable of a given optimization problem is referred to as \emph{implicit} whenever
it is used in order to model feasibility conditions but does not appear in the
objective function. On the one hand, such variables seem to be of less interest for
the actual purpose of optimization but, on the other hand, offer the possibility to transfer
difficult constraints into feasibility conditions which can be handled easier. 

In this paper, we are concerned with the model problem
\begin{equation}\label{eq:implicit_problem}\tag{P}
	\min\limits_z\{f(z)\,|\,z\in M\cap\dom K\}
\end{equation}
where $f\colon\R^n\to\R$ is, for simplicity, a continuously differentiable function,
$M\subset\R^n$ is a closed set modeling standard constraints,
$K\colon\R^n\tto\R^m$ is a set-valued mapping with a closed graph,
and $\dom K$ denotes the domain of $K$, i.e., the set of points
such that the associated images of $K$ are nonempty.
Often, model \eqref{eq:implicit_problem} is difficult to deal with due to the presence
of the implicit constraint $z\in\dom K$ as, typically, $K$ is of challenging variational
structure. Noting that this constraint demands the existence of some $\lambda\in K(z)$
motivates the investigation of the explicit problem
\begin{equation}\label{eq:explicit_problem}\tag{Q}
	\min\limits_{z,\lambda}\{f(z)\,|\,z\in M,\,(z,\lambda)\in\gph K\}
\end{equation}
which possesses the extra variable $\lambda$.
Above, $\gph K$ denotes the graph of $K$.
Following \cite{BenkoMehlitz2021}, where a more complicated model is considered,
we refer to $\lambda$ as an \emph{implicit} variable in \eqref{eq:implicit_problem} as
it is implicitly hidden in the constraint system.
In \eqref{eq:explicit_problem}, $\lambda$ is \emph{explicitly} used as a variable.

Let us illustrate the seemingly beneficial use of implicit variables in 
the modeling of optimization problems.
We consider cardinality-constrained optimization problems.
Note that the cardinality constraint
\begin{equation}\label{eq:cardinality_constraint}
	\norm{z}_0\leq\kappa,
\end{equation}
where $\norm{\cdot}_0\colon\R^n\to\{0,1,\ldots,n\}$, the so-called $\ell_0$-quasi-norm, 
counts the nonzero entries of the argument vector
and $\kappa\in\{1,\ldots,n-1\}$ is a given constant,
is, in certain sense, equivalent to the complementarity-type system
\begin{equation}\label{eq:ref_complementarity}
\begin{aligned}
	\mathtt e^\top \lambda&\geq n-\kappa\\
	0\leq \lambda&\leq\mathtt e\\
	z\bullet \lambda&=0,
\end{aligned}
\end{equation}
see \cite{BurdakovKanzowSchwartz2016} for a detailed study.
Above, $\mathtt e\in\R^n$ is the all-ones vector and $z\bullet \lambda$ denotes the
entrywise product of the vectors $z,\lambda\in\R^n$.
In numerical practice, it seems to be beneficial to use the reformulation
\eqref{eq:ref_complementarity} of the 
variationally more challenging original constraint \eqref{eq:cardinality_constraint} since complementarity-type
systems in the constraints of an optimization problem
can be treated with already available penalty, multiplier penalty, relaxation, active set, SQP-, or Newton-type methods,
see e.g.\ \cite{FletcherLeyfferRalphScholtes2006,GuoDeng2021,HarderMehlitzWachsmuth2021,HoheiselKanzowSchwartz2013,HuangYangZhu2006,
IzmailovSolodov2008,JudiceSheraliRibeiroFaustino2007,LeyfferLopezNocedal2006} and the references therein.
This particular situation is covered by our modeling approach when using $K\colon\R^n\tto\R^n$ given by
\[
	\forall z\in\R^n\colon\quad
	K(z):=\{\lambda\in\R^n\,|\,\mathtt e^\top\lambda\geq n-\kappa,\,0\leq\lambda\leq\mathtt e,\,z\bullet\lambda=0\},
\]
see \cref{sec:example_CCP}.
Hence, the approach promoted in \cite{BurdakovKanzowSchwartz2016} exploits $\lambda$ as an \emph{explicit} variable
in order to avoid dealing with the challenging cardinality constraint \eqref{eq:cardinality_constraint} directly.
In the recent paper \cite{BenkoMehlitz2021}, it is shown theoretically and in terms of a more complex model problem
that this approach comes along with several drawbacks.
Potentially, the reformulated problem possesses more stationary and, thus, locally
optimal points. Furthermore, one may need more restrictive constraint qualifications for
the reformulated problem than for the original one in order to derive comparable
necessary optimality conditions. Finally, one has to deal with more variables in the reformulated problem.
Let us mention that the phenomenon of implicit variables has some major impact in the context
of bilevel and evaluated multiobjective optimization as well. In these settings, 
lower-level Lagrange multipliers and scalarization parameters play the role of implicit variables, respectively,
see e.g.\ \cite{AdamHenrionOutrata2018,DempeDutta2012,DempeMehlitz2020,DempeMehlitz2025}
and \cref{sec:example_BPP} below.

The purpose of this paper is twofold.
First, as already mentioned before, in \cite{BenkoMehlitz2021},
the authors illustrated the issues caused by treating implicit variables as explicit ones
in terms of a rather complex model problem.
The present paper aims to recover these findings for the much simpler and easily accessible
model problem \eqref{eq:implicit_problem}.
Furthermore, in \cite{BenkoMehlitz2021}, when comparing the implicit and explicit problem
with respect to (w.r.t.) stationary points, the authors merely investigated so-called
Mordukhovich-stationary points which are based on the limiting tools of variational analysis,
see e.g.\ \cite{Mordukhovich2018}.
Here, we extend these findings to so-called Bouligand-stationary and strongly stationary points.
Therefore, our main workhorse is \cite[Theorem~3.1]{BenkoMehlitz2022} which explains quite nicely
the relations between tangents and normals to the domain and the graph of a set-valued mapping.
Second, we aim to demonstrate by means of a simple computational experiment
that treating implicit variables as explicit ones may also lead to numerical issues. 
Therefore, we investigate portfolio optimization problems with cardinality constraints
of type \eqref{eq:cardinality_constraint} and their associated reformulation based on the
complementarity-type system \eqref{eq:ref_complementarity}.
For a fair comparison, which is based on a benchmark collection of example problems,
we rely on the augmented Lagrangian framework from \cite{JiaKanzowMehlitzWachsmuth2023}
which nicely applies to the original model and its reformulation.
Let us mention that \cite{BenkoMehlitz2021} does not present any computational results.

Let us comment on already available literature which points out that treating implicit variables as explicit
ones might be numerically disadvantageous. In \cite{Mehlitz2020c}, the author compares different
solution methods for so-called or-constrained optimization problems which, for the price of additional
slack variables and constraints, can be reformulated as complementarity- or switching-constrained 
optimization problems. Exemplary, all of these problems can be treated with the aid of similar relaxation
methods, see e.g.\ \cite{HoheiselKanzowSchwartz2013,KanzowMehlitzSteck2021} as well. 
However, by means of numerical experiments, it has been shown that, apart from artificial
situations, a direct treatment of the original model gives the best results.
In \cite{ZemkohoZhou2021}, a numerical comparison of the so-called value function and 
Karush--Kuhn--Tucker reformulation of bilevel optimization problems is provided. While the
latter comes along with additional variables in form of lower-level Lagrange multipliers,
the former possesses the same number of variables as the original problem. 
Both reformulations can be
solved by applying a semismooth Newton-type method to the resulting stationarity systems. 
It is reported in \cite{ZemkohoZhou2021} that the value function reformulation outperforms the
Karush--Kuhn--Tucker reformulation. This is not only caused by the higher number of variables and
constraints in the Karush--Kuhn--Tucker reformulation
but also results from the fact that derivatives of higher order 
are necessary in order to tackle the latter.
The numerical experiments in the present paper similarly will underline that treating implicit variables
as explicit ones worsens the computational performance when comparable solution algorithms are exploited
to tackle the original problem and its reformulation.

The remainder of the paper is structured as follows.
In \cref{sec:notation}, we comment on the notation which is exploited in this paper.
Furthermore, we recall some essential concepts from variational analysis
and discuss associated preliminary result.
\cref{sec:examples} aims to motivate the consideration of the model problem \eqref{eq:implicit_problem}
and its reformulation \eqref{eq:explicit_problem} by means of three different example classes,
namely cardinality-constrained, vanishing-constrained, and bilevel optimization problems.
The theoretical analysis of the relationship between \eqref{eq:implicit_problem} and \eqref{eq:explicit_problem}
is carried out in \cref{sec:analysis}.
\cref{sec:minimizers,sec:relations_stationary_points} focus on a comparison of both
problems in terms of minimizers and stationarity points as well as constraint qualifications, respectively.
A compact summary is presented in \cref{sec:summary}.
\cref{sec:comparison} interrelates these findings and the ones from \cite{BenkoMehlitz2021}
which were obtained in terms of a more complex model problem which is covered by \eqref{eq:implicit_problem}.
In \cref{sec:consequences_examples}, we revisit the example problems from \cref{sec:examples}
and illustrate some of our theoretical findings from \cref{sec:analysis}.
The aforementioned numerical experiments in terms of cardinality-constrained portfolio optimization problems
are stated in \cref{sec:numerics}.
Some concluding remarks close the paper in \cref{sec:conclusions}.

\section{Notation and preliminaries}\label{sec:notation}

In this section, we comment on the notation which is used in this paper.
Furthermore, we recall some standard constructions from variational analysis
and present associated preliminary results which address the calculus
of the variational objects of our interest.

\subsection{Basic notation}\label{sec:basic_notation}

Let $\N$, $\R$, $\R_+$, and $\R_-$ denote the positive integers
as well as the real, nonnegative real, and nonpositive real numbers, respectively.
Given $z,w\in\R^n$, $z\bullet w\in\R^n$ is the componentwise product of $z$ and $w$.
Furthermore, $\norm{z}$ denotes the Euclidean norm of $z$.
For a nonempty set $\Omega\subset\R^n$, $\dist(z,\Omega):=\inf_{z'\in\Omega}\norm{z'-z}$
is the distance of $z$ to $\Omega$, 
and $\dist(\cdot,\Omega)\colon\R^n\to\R$ denotes the associated distance function of $\Omega$.
Given a differentiable mapping $\upsilon\colon\R^n\to\R^m$,
$\upsilon'(z)\in\R^{m\times n}$ represents the Jacobian of $\upsilon$ at $z$.
Whenever $m:=1$ is valid, we exploit $\nabla\upsilon(z):=\upsilon'(z)^\top$
to denote the gradient of $\upsilon$ at $z$.
Partial derivatives w.r.t.\ certain variables are denotes in canonical way.

\subsection{Variational analysis}\label{sec:variational_analysis}

In this subsection,
we review some fundamental concepts from variational analysis
which can be found in standard text books like \cite{Mordukhovich2018,RockafellarWets1998}.
Furthermore, we recall some Lipschitzian properties of set-valued mappings from \cite{Benko2021}.

For a set $\Omega\subset\R^n$ and some point $\bar z\in\Omega$ such that $\Omega$ is closed
locally around $\bar z$,
\[
	T_\Omega(\bar z)
	:=
	\left\{
		w\in\R^n\,\middle|\,
		\begin{aligned}
			&\exists\{t_k\}_{k\in\N}\subset(0,\infty)\,\exists\{w_k\}_{k\in\N}\subset\R^n\colon
			\\
			&\quad t_k\downarrow 0,\,w_k\to w,\,\bar z+t_kw_k\in\Omega\,\forall k\in\N
		\end{aligned}
	\right\}
\]
denotes the standard (Bouligand) tangent cone to $\Omega$ at $\bar z$.
Furthermore, we are going to exploit
the regular (Fr\'{e}chet) normal cone to $\Omega$ at $\bar z$, given by
\[
	\widehat N_\Omega(\bar z)
	:=
	\{
		\zeta\in\R^n\,|\,\forall z\in\Omega\colon\,\zeta^\top(z-\bar z)\leq\oo(\norm{z-\bar z})
	\},
\]
and the limiting (Mordukhovich) normal cone to $\Omega$ at $\bar z$, given by
\[
	N_\Omega(\bar z)
	:=
	\left\{
		\zeta\in\R^n\,\middle|\,
		\begin{aligned}
			&\exists\{z_k\}_{k\in\N}\subset\Omega\,\exists\{\zeta_k\}_{k\in\N}\subset\R^n\colon
			\\
			&\quad z_k\to\bar z,\,\zeta_k\to\zeta,\,\zeta_k\in\widehat N_\Omega(z_k)\,\forall k\in\N
		\end{aligned}
	\right\}.
\]
It is well known that the regular normal cone is the polar of the tangent cone, i.e.,
\[
	\widehat N_\Omega(\bar z)
	=
	\{\zeta\in\R^n\,|\,\forall w\in T_\Omega(\bar z)\colon\,\zeta^\top w\leq 0\}.
\]
By definition, we have $\widehat N_\Omega(\bar z)\subset N_\Omega(\bar z)$ in general,
and if equality holds, $\Omega$ is called regular at $\bar z$.
Whenever $\Omega$ is convex, it is regular at each of its points, 
and the regular and limiting normal cone both coincide
with the normal cone of convex analysis.

For a lower semicontinuous function $\upsilon\colon\R^n\to\R$ and some point $\bar z\in\R^n$,
we refer to
\[
	\widehat\partial\upsilon(\bar z)
	:=
	\{\zeta\in\R^n\,|\,(\zeta,-1)\in\widehat N_{\epi\upsilon}(\bar z,\upsilon(\bar z))\}
\]
as the regular subdifferential of $\upsilon$ at $\bar z$.
Here, $\epi\upsilon:=\{(z,\alpha)\in\R^n\times\R\,|\,\alpha\geq\upsilon(z)\}$
is the epigraph of $\upsilon$.

For a set-valued mapping $\Upsilon\colon\R^n\tto\R^m$,
we exploit $\gph\Upsilon:=\{(x,y)\in\R^n\times\R^m\,|\,y\in\Upsilon(x)\}$ 
and $\dom\Upsilon:=\{x\in\R^n\,|\,\Upsilon(x)\neq\emptyset\}$
to denote the graph and the domain of $\Upsilon$.
Furthermore, $\Upsilon^{-1}\colon\R^m\tto\R^n$ given by 
$\gph\Upsilon^{-1}:=\{(y,x)\in\R^m\times\R^n\,|\,(x,y)\in\gph\Upsilon\}$
is the inverse of $\Upsilon$.
We call $\Upsilon$ a polyhedral mapping if its graph can be represented as the union
of finitely many convex polyhedral sets.

Given some point $\bar x\in\dom\Upsilon$, $\Upsilon$ is called
inner semicompact w.r.t.\ $\dom\Upsilon$ at $\bar x$ whenever for each
sequence $\{x_k\}_{k\in\N}\subset\dom\Upsilon$ such that $x_k\to\bar x$,
there is a sequence $\{y_k\}_{k\in\N}\subset\R^m$ with a bounded subsequence
such that $y_k\in\Upsilon(x_k)$ is valid for each $k\in\N$.
We note that each mapping that possesses uniformly bounded images around
the reference point in its domain is inner semicompact w.r.t.\ its domain there.
It is important to note that whenever $\Upsilon$ possesses a closed graph 
while being inner semicompact w.r.t.\ its domain at $\bar x\in\dom\Upsilon$,
then $\dom\Upsilon$ is closed locally around $\bar x$,
see \cite[Lemma~2.1]{Benko2021}.
In this paper, this observation will be used frequently.

The mapping $\Upsilon$ is said to be
inner calm* in the fuzzy sense w.r.t.\ $\dom\Upsilon$ at $\bar x$
whenever for each $d\in T_{\dom\Upsilon}(\bar x)$ with $\norm{d}=1$,
there exist a constant $\kappa_d>0$, some point $\bar y\in\R^m$, as well as sequences
$\{t_k\}_{k\in\N}\subset(0,\infty)$, $\{d_k\}_{k\in\N}\subset\R^n$,
and $\{y_k\}_{k\in\N}\subset\R^m$
such that $t_k\downarrow 0$, $d_k\to d$, and $y_k\to\bar y$
as well as
$y_k\in\Upsilon(\bar x+t_kd_k)$ and $\norm{y_k-\bar y}\leq\kappa_dt_k\,\norm{d_k}$
for each $k\in\N$ are valid.
Let us note that this property, originating from \cite{Benko2021},
serves as a quantitative version of inner semicompactness of a mapping.

At some point $(\bar x,\bar y)\in\gph\Upsilon$, 
$\Upsilon$ is called metrically subregular if there exist a constant $\kappa>0$
and a neighborhood $U\subset\R^n$ of $\bar x$ such that
\[
	\forall x\in U\colon\quad 
	\dist(x,\Upsilon^{-1}(\bar y))
	\leq
	\kappa\,\dist(\bar y,\Upsilon(x)).
\]
From \cite[Theorem~3.4]{Benko2021} we know that polyhedral mappings are 
inner calm* in the fuzzy sense w.r.t.\ their domain at each point of their domain
and metrically subregular at each point of their graph.

For $(\bar x,\bar y)\in\gph\Upsilon$, the set-valued mapping $D\Upsilon(\bar x,\bar y)\colon\R^n\tto\R^m$
given by
\[
	\gph D\Upsilon(\bar x,\bar y):=T_{\gph\Upsilon}(\bar x,\bar y)
\]
is referred to as the graphical derivative of $\Upsilon$ at $(\bar x,\bar y)$.
Furthermore, the set-valued mappings $\widehat D^*\Upsilon(\bar x,\bar y)\colon\R^m\tto\R^n$ and
$D^*\Upsilon(\bar x,\bar y)\colon\R^m\tto\R^n$ defined via
\begin{align*}
	\forall y^*\in\R^m\colon\quad
	\widehat D^*\Upsilon(\bar x,\bar y)(y^*)
	&:=
	\{x^*\in\R^n\,|\,(x^*,-y^*)\in \widehat N_{\gph\Upsilon}(\bar x,\bar y)\},
	\\
	D^*\Upsilon(\bar x,\bar y)(y^*)
	&:=
	\{x^*\in\R^n\,|\,(x^*,-y^*)\in N_{\gph\Upsilon}(\bar x,\bar y)\}
\end{align*}
are called the regular and limiting coderivative of $\Upsilon$ at $(\bar x,\bar y)$, respectively.

\subsection{Preliminary results}\label{sec:preliminaries}

In the course of the paper,
the following lemmas are of essential importance.

To start, we aim to clarify how tangent and normal cones to the domain of a set-valued
mapping can be computed or, at least, estimated in terms of its graphical derivative
and coderivatives.
This result is taken from \cite[Theorem~3.1]{BenkoMehlitz2022}.

\begin{lemma}\label{lem:dom_graph_calculus}
	Let $K\colon\R^n\tto\R^m$ be a set-valued mapping with a closed graph,
	and fix $\bar z\in\dom K$ where $K$ is inner semicompact w.r.t.\ $\dom K$.
	Then the following assertions hold.
	\begin{enumerate}
		\item We always have
		\[
			T_{\dom K}(\bar z)
			\supset
			\bigcup\limits_{\bar\lambda\in K(\bar z)}\dom D K(\bar z,\bar\lambda),
		\]
		and the converse inclusion holds whenever $K$ is inner calm* in the fuzzy sense
		at $\bar z$ w.r.t.\ $\dom K$.
		\item We always have
		\[	
			\widehat N_{\dom K}(\bar z)
			\subset
			\bigcap\limits_{\bar\lambda\in K(\bar z)}\widehat D^* K(\bar z,\bar\lambda)(0),
		\]
		and the converse inclusion holds whenever $K$ is inner calm* in the fuzzy sense
		at $\bar z$ w.r.t.\ $\dom K$.
		\item We have
		\[
			N_{\dom K}(\bar z)
			\subset
			\bigcup\limits_{\bar\lambda\in K(\bar z)} D^* K(\bar z,\bar\lambda)(0).
		\]
	\end{enumerate}
\end{lemma}

Let us recall that the inner semicompactness of $K$ at $\bar z$ w.r.t.\ $\dom K$,
which is an assumption in \cref{lem:dom_graph_calculus}, guarantees that $\dom K$
is locally closed around $\bar z$.

	Below, we present an intersection rule.
	The one for tangents and limiting normals is taken from \cite[Section~5.1.2]{BenkoMehlitz2022}
	while the one for regular normals can be found in \cite[Theorem~3.3]{HuongAnXu2024}.

\begin{lemma}\label{lem:intersection_rule}
	Let $\Omega_1,\Omega_2\subset\R^n$ be sets and fix $\bar z\in\Omega_1\cap\Omega_2$
	such that $\Omega_1$ and $\Omega_2$ are locally closed around $\bar z$.
	Then the following assertions hold.
	\begin{enumerate}
		\item We always have
			\[
				T_{\Omega_1\cap\Omega_2}(\bar z)
				\subset
				T_{\Omega_1}(\bar z)\cap T_{\Omega_2}(\bar z),
			\]
			and the converse inclusion holds whenever the mapping
			$z\mapsto(z,z)-\Omega_1\times\Omega_2$ is metrically subregular at
			$(\bar z,(0,0))$.
		\item We always have
			\[
				\widehat N_{\Omega_1\cap\Omega_2}(\bar z)
				\supset
				\widehat N_{\Omega_1}(\bar z)+\widehat N_{\Omega_2}(\bar z),
			\]
			and the converse inclusion holds if and only if there is a constant
			$\kappa>0$ such that
			\begin{equation}\label{eq:abstract_CQ_for_S_stat}
				\widehat\partial\dist(\cdot,\Omega_1\cap\Omega_2)(\bar x)
				\subset
				\kappa\bigl(
					\widehat\partial\dist(\cdot,\Omega_1)(\bar x)
					+
					\widehat\partial\dist(\cdot,\Omega_2)(\bar x)
				\bigr).
			\end{equation}
		\item If the mapping $z\mapsto(z,z)-\Omega_1\times\Omega_2$ 
			is metrically subregular at $(\bar z,(0,0))$, then we have
			\[
				N_{\Omega_1\cap\Omega_2}(\bar z)
				\subset
				N_{\Omega_1}(\bar z) + N_{\Omega_2}(\bar z).
			\]
	\end{enumerate}
\end{lemma}

We note that the subregularity requirement in \cref{lem:intersection_rule} boils down
to the existence of a neighborhood $U\subset\R^n$ of the reference point $\bar z$
and a constant $\kappa>0$ such that
\[
	\forall z\in U\colon\quad
	\dist(z,\Omega_1\cap\Omega_2)
	\leq
	\kappa(\dist(z,\Omega_1)+\dist(z,\Omega_2)),
\]
the latter being well known as subtransversality of $\Omega_1$ and $\Omega_2$ at $\bar z$,
and it is common knowledge that the calculus rules in \cref{lem:intersection_rule},
which address tangents and limiting normals, hold in the presence of this condition, 
see e.g.\ \cite[Section~4]{KrugerLukeThao2018} and the references therein.
Furthermore, it should be noted that subtransversality is not enough to guarantee
validity of the inclusion $\subset$ in the intersection rule for regular normals
as the following simple example indicates.

\begin{example}\label{ex:intersection_rule_regular_normals}
	Consider the sets
	\[
		\Omega_1:=\{z\in\R^2_+\,|\,z_1z_2=0\},
		\qquad
		\Omega_2:=\{z\in\R^2\,|\,z_1=z_2\}.
	\]
	On the one hand, we have $\Omega_1\cap\Omega_2=\{0\}$ yielding
	\[
		\widehat N_{\Omega_1\cap\Omega_2}(0)=\R^2.
	\]
	On the other hand, one can easily check that
	\[
		\widehat N_{\Omega_1}(0)=\R^2_-,
		\qquad
		\widehat N_{\Omega_2}(0)=\{\zeta\in\R^2\,|\,\zeta_1+\zeta_2=0\},
	\]
	and this gives
	\[
		\widehat N_{\Omega_1}(0)+\widehat N_{\Omega_2}(0)
		=
		\{\zeta\in\R^2\,|\,\zeta_1+\zeta_2\leq 0\}.
	\]
	Furthermore, one should note that $\Omega_1$ and $\Omega_2$ are polyhedral sets.
	Hence, the set-valued mapping $z\mapsto(z,z)-\Omega_1\times\Omega_2$ is polyhedral and, thus,
	metrically subregular at each point of its graph,
	yielding subtransversality of $\Omega_1$ and $\Omega_2$ at $0$.
\end{example}

	Following e.g.\ \cite[Section~5.1]{BenkoMehlitz2022},
	a sufficient condition for the metric subregularity assumption 
	in \cref{lem:intersection_rule} is given by
	\begin{equation}\label{eq:Mord_crit_intersection}
		N_{\Omega_1}(\bar z)\cap\bigl(-N_{\Omega_2}(\bar z)\bigr)=\{0\}.
	\end{equation}
	To illustrate this condition, 
	let us assume that $\Omega_i=\{z\in\R^n\,|\,\upsilon_i(z)\in C_i\}$
	holds for $i=1,2$, where $\upsilon_i\colon\R^n\to\R^{m_i}$ is continuously differentiable,
	$C_i\subset\R^{m_i}$ is closed, and $m_i\in\N$.
	Then \eqref{eq:Mord_crit_intersection} is implied by
	\[
		\upsilon'_1(\bar z)^\top\eta_1+\upsilon'_2(\bar z)^\top \eta_2 = 0,
		\,
		\eta_i\in N_{C_i}(\upsilon_i(\bar z))\,(i=1,2)
		\quad\Longrightarrow\quad
		\eta_i=0\,(i=1,2),
	\]
	see e.g.\ \cite[Theorem~6.14]{RockafellarWets1998},
	which can be interpreted as a generalized version of the
	Mangasarian--Fromovitz constraint qualification.
	In the literature,
	this type of condition is often referred to as
	no nonzero abnormal multiplier constraint qualification (NNAMCQ for short).
	Furthermore, in this exemplary setting,
	\eqref{eq:abstract_CQ_for_S_stat} holds whenever
	\[
		\begin{bmatrix}
			\upsilon_1'(\bar z)
			\\
			\upsilon_2'(\bar z)
		\end{bmatrix}
		\R^n
		+
		L
		=
		\R^{m_1+m_2}
	\]
	is valid for some subspace $L\subset\R^{m_1+m_2}$ which satisfies
	\[
		T_{C_1\times C_2}(\upsilon_1(\bar z),\upsilon_2(\bar z))
		+
		L
		\subset
		T_{C_1\times C_2}(\upsilon_1(\bar z),\upsilon_2(\bar z)),
	\]
	see \cite[Theorem~4]{GfrererOutrata2016a}.
	Clearly, the latter holds if the matrix
	\[
		\begin{bmatrix}
			\upsilon_1'(\bar z)
			\\
			\upsilon_2'(\bar z)
		\end{bmatrix}
		\in\R^{(m_1+m_2)\times n}
	\]
	possesses full row rank $m_1+m_2$ as this allows for $L$ being trivial.
	More enhanced choices of $L$ lead to conditions
	related to the linear independence constraint qualification (LICQ for short),
	see e.g.\ \cite{Mehlitz2020b} for an illustration in terms of disjunctive optimization.

\section{Motivating examples}\label{sec:examples}

In this section,
we present three classes of practically relevant optimization problems \eqref{eq:implicit_problem} 
which can be modeled with the aid of implicit variables.

\subsection{Cardinality-constrained optimization}\label{sec:example_CCP}

For $n\in\N$ and $\kappa\in\{1,\ldots,n-1\}$, we define
\[
	\R_{\leq\kappa}^n:=\{z\in\R^n\,|\,\norm{z}_0\leq\kappa\}.
\]
Recall that $\norm{\cdot}_0\colon\R^n\to\{0,1,\ldots,n\}$ is the function which counts the
nonzero entries of the input vector. 
We refer to
\begin{equation}\label{eq:CCP}\tag{P$_\textup{cc}$}
	\min\limits_z\{f(z)\,|\,z\in M,\,z\in\R_{\leq\kappa}^n\}
\end{equation}
as a cardinality-constrained optimization problem. 
Cardinality constraints typically appear in optimization problems related to 
applications in data science like data compression or sparse recovery.
In \cite{BurdakovKanzowSchwartz2016}, the authors suggest to reformulate the
cardinality constraint in \eqref{eq:CCP} with the aid of additional variables.
More precisely, they investigate the surrogate model
\begin{equation}\label{eq:CCPref}\tag{Q$_\textup{cc}$}
	\min\limits_{z,\lambda}
	\{f(z)\,|\,
		z\in M,\,\mathtt e^\top \lambda\geq n-\kappa,\,0\leq\lambda\leq\mathtt e,\,
		z\bullet\lambda=0
	\},
\end{equation}
whose constraints possess complementarity-type structure and, thus, allow for
the derivation of optimality conditions and solution algorithms.
The reformulated problem \eqref{eq:CCPref} has been used e.g.\ in
\cite{BucherSchwartz2018,CervinkaKanzowSchwartz2016,KanzowRaharjaSchwartz2021a}
for the derivation of optimality conditions, constraint qualifications, and solution algorithms
for \eqref{eq:CCP}.

Introducing $K_\textup{cc}\colon\R^n\tto\R^n$ by means of
\begin{equation}\label{eq:def_Kcc}
	\forall z\in\R^n\colon\quad
	K_\textup{cc}(z)
	:=
	\{\lambda\in\R^n\,|\,\mathtt e^\top \lambda\geq n-\kappa,\,0\leq\lambda\leq\mathtt e,\,
		z\bullet\lambda=0\},
\end{equation}
\cite[Theorem~3.1]{BurdakovKanzowSchwartz2016} guarantees
\[
	\dom K_\textup{cc}=\R^n_{\leq\kappa},
\]
so \eqref{eq:CCP} fits perfectly into our setting.

\subsection{Optimization problems with vanishing constraints}\label{sec:example_MPVC}

Given continuously differentiable functions $G,H\colon\R^n\to\R^m$, 
an optimization problem of the form
\begin{equation}\label{eq:MPVC}\tag{P$_\textup{vc}$}
	\min\limits_z\{f(z)\,|\,z\in M,\,H(z)\geq 0,\,G(z)\bullet H(z)\leq 0\}
\end{equation}
is referred to as a mathematical problem with vanishing constraints.
Exemplary, these models are used to optimize truss structures, 
see \cite{AchtzigerKanzow2008} for details.
Further applications, theory, and numerical methods associated with vanishing-constrained
optimization problems can be found e.g.\ in 
\cite{AchtzigerKanzow2008,Hoheisel2009,HoheiselPablosPooladianSchwartzSteverango2020,HoheiselKanzow2007}.
Due to the fact that complementarity-constrained problems are much more frequently studied
in the literature,
it has been mentioned in \cite[Section~2]{AchtzigerKanzow2008} that
\eqref{eq:MPVC} can be accessed via the reformulation
\begin{equation}\label{eq:MPVCref}\tag{Q$_\textup{vc}$}
	\min\limits_{z,\lambda}
	\left\{f(z)\,\middle|\,
		\begin{aligned}
			&z\in M,\,G(z)-\lambda\leq 0,
			\\
			&H(z)\geq 0,\,\lambda\geq 0,\,H(z)\bullet\lambda=0
		\end{aligned}
	\right\}
\end{equation}
which is indeed a mathematical problem with complementarity constraints.
Exploiting $K_\textup{vc}\colon\R^n\tto\R^m$ given by
\begin{equation}\label{eq:def_Kvc}
	\forall z\in\R^n\colon\quad
	K_\textup{vc}(z)
	:=\{\lambda\in\R^m\,|\,G(z)-\lambda\leq 0,\,H(z)\geq 0,\,
		\lambda\geq 0,\,H(z)\bullet\lambda=0
		\},
\end{equation}
one can easily check that
\[
	\dom K_\textup{vc} 
	=
	\{z\in\R^n\,|\,H(z)\geq 0,\,G(z)\bullet H(z)\leq 0\}.
\]
Indeed, $z\in\dom K_\textup{vc}$ yields the existence of $\lambda\in\R^m$ such that
$\lambda\geq\max(G(z),0)$, $H(z)\geq 0$, and $H(z)\bullet\lambda=0$.
For each $i\in\{1,\ldots,m\}$, $\lambda_i=0$ means $G_i(z)\leq 0$ and, hence,
$G_i(z)H_i(z)\leq 0$ while $\lambda_i>0$ gives $H_i(z)=0$ and, thus,
$G_i(z)H_i(z) = 0$.
Conversely, given $z\in\R^n$ such that $H(z)\geq 0$ and $G(z)\bullet H(z)\leq 0$,
one can set $\lambda:=\max(G(z),0)$ in order to find $\lambda\in K_\textup{vc}(z)$,
i.e., $z\in\dom K_\textup{vc}$, see \cite[Lemma~1]{AchtzigerKanzow2008} as well.

\subsection{Bilevel optimization}\label{sec:example_BPP}

Fix $n_1,n_2\in\N$ such that $n_1+n_2=n$.
Throughout the section, we exploit $\R^n=\R^{n_1}\times\R^{n_2}$
and split $z\in\R^n$ into $(z_1,z_2)\in\R^{n_1}\times\R^{n_2}$.
Given continuously differentiable functions $g\colon\R^{n}\to\R$ and $G\colon\R^n\to\R^m$
such that, for each $z_1\in\R^{n_1}$, $g(z_1,\cdot)$ is convex 
while $G(z_1,\cdot)$ is componentwise convex, 
we investigate the parametric optimization problem
\begin{equation}\label{eq:lower_level}\tag{Par$(z_1)$}
	\min\limits_{z_2}\{g(z)\,|\,G(z)\leq 0\}.
\end{equation}
Let us define the associated solution mapping $\Psi\colon\R^{n_1}\tto\R^{n_2}$ by means of
\begin{align*}
	\forall z_1\in\R^{n_1}\colon\quad
	\Psi(z_1):=\argmin\limits_{z_2}\{g(z)\,|\,G(z)\leq 0\}.
\end{align*}
We are now considering the (optimistic) bilevel optimization problem
\begin{equation}\label{eq:BOP}\tag{P$_\textup{bop}$}
	\min\limits_{z}\{f(z)\,|\,z\in M,\,z\in\gph\Psi\}
\end{equation}
which possesses numerous applications in data, management, and natural sciences
as well as economics. A detailed introduction to bilevel optimization can be found
in the monographs \cite{Bard1998,Dempe2002,DempeKalashnikovPerezValdesKalashnykova2015}.

A popular approach to transfer \eqref{eq:BOP} into an accessible single-level optimization
problem exploits the Karush--Kuhn--Tucker conditions of \eqref{eq:lower_level}.
In order to apply it, let us assume that the convex optimization problem \eqref{eq:lower_level}
satisfies Slater's constraint qualification for each 
$z_1\in\dom\Psi$.
Then, for any such $z_1$, relation $z\in\gph\Psi$ is equivalent 
to the existence of a multiplier
$\lambda\in\R^m$ such that
\[
	\nabla_{z_2}g(z)+G'_{z_2}(z)^\top\lambda=0,\,
	G(z)\leq 0,\,\lambda\geq 0,\,G(z)\bullet\lambda=0.
\]
This motivates the investigation of
\begin{equation}\label{eq:BOPref}\tag{Q$_\textup{bop}$}
	\min\limits_{z,\lambda}\left\{f(z)\,\middle|\,
		\begin{aligned}
		&z\in M,\,\nabla_{z_2}g(z)+G'_{z_2}(z)^\top\lambda=0,
		\\
		&G(z)\leq 0,\,\lambda\geq 0,\,G(z)\bullet\lambda=0
		\end{aligned}
		\right\}
\end{equation}
which, again, is an optimization problem with complementarity constraints.
We would like to emphasize that this so-called 
Karush--Kuhn--Tucker reformulation of \eqref{eq:BOP}
is used frequently in the literature to tackle the latter,
see e.g.\ \cite{DempeZemkoho2012,KimLeyfferMunson2020} and the references therein,
as complementarity-constrained optimization problems are, 
as already mentioned in \cref{sec:introduction}, 
despite being notoriously difficult, 
well-understood from a theoretical and numerical point of view.

Defining $K_\textup{bop}\colon\R^{n}\tto\R^m$ by means of
\begin{equation}\label{eq:def_Kbop}
	\forall z\in\R^{n}\colon\quad
	K_\textup{bop}(z)
	:=
	\left\{\lambda\in\R^m\,\middle|\,
		\begin{aligned}
		&\nabla_{z_2}g(z)+G'_{z_2}(z)^\top\lambda=0,
		\\
		&G(z)\leq 0,\,\lambda\geq 0,\,G(z)\bullet\lambda=0
		\end{aligned}
	\right\},
\end{equation}
our modeling approach covers the above bilevel optimization problem \eqref{eq:BOP}
as we already outlined above that $\dom K_\textup{bop}=\gph\Psi$ is valid.

For the purpose of completeness,
let us mention that the use of extra variables is very popular in hierarchical optimization
in order to transfer the original model into a more tractable one.
Exemplary, similar as above, one can replace the implicit constraint $z\in\gph\Psi$
in \eqref{eq:BOP} by a relation demanding strong duality (w.r.t.\ a suitable duality concept)
for \eqref{eq:lower_level}. Then the variables of the dual problem associated with 
\eqref{eq:lower_level} play the role of implicit variables,
see \cite{DempeMehlitz2025} for a recent survey.
When considering a bilevel optimization problem whose underlying parametric optimization
problem possesses a vector-valued objective function,
it is possible to exploit a scalarization approach in order to transfer it into a
bilevel optimization problem of type \eqref{eq:BOP}.
However, the scalarization parameters then play the role of implicit variables already,
see \cite{DempeMehlitz2020}.

\section{Analysis of the model problem}\label{sec:analysis}

This section is devoted to a theoretical analysis which compares the
model problem \eqref{eq:implicit_problem} and its reformulation \eqref{eq:explicit_problem}
w.r.t.\ their minimizers, stationary points, and (selected) constraint qualifications.
These results are collected in \cref{sec:minimizers,sec:relations_stationary_points},
and a brief summary is presented in \cref{sec:summary}.
In \cref{sec:comparison}, we embed our findings into the ones from \cite{BenkoMehlitz2021}.

\subsection{Relations of minimizers}\label{sec:minimizers}

The following result, 
which discusses the relationship between global minimizers 
of \eqref{eq:implicit_problem} and \eqref{eq:explicit_problem}, 
is a simple consequence of the definition of the domain of a set-valued mapping 
which is why we omit the associated straightforward proof.
\begin{theorem}\label{thm:global_minimizers}
	\phantom{firstline}
	\begin{enumerate}
		\item If $\bar z\in\R^n$ is a global minimizer of \eqref{eq:implicit_problem},
			then $(\bar z,\bar\lambda)$ is a global minimizer of \eqref{eq:explicit_problem}
			for each $\bar\lambda\in K(\bar z)$.
		\item If $(\bar z,\bar\lambda)\in\R^n\times\R^m$ is a global minimizer 
			of \eqref{eq:explicit_problem}, then $\bar z$ is a global minimizer 
			of \eqref{eq:implicit_problem}.
	\end{enumerate}
\end{theorem}

Our next result addresses the relationship between the local minimizers 
of \eqref{eq:implicit_problem} and \eqref{eq:explicit_problem},
and is more delicate than \cref{thm:global_minimizers}.
\begin{theorem}\label{thm:local_minimizers}
	\phantom{firstline}
	\begin{enumerate}
		\item If $\bar z\in\R^n$ is a local minimizer of \eqref{eq:implicit_problem},
			then $(\bar z,\bar\lambda)$ is a local minimizer of \eqref{eq:explicit_problem}
			for each $\bar\lambda\in K(\bar z)$.
		\item If $(\bar z,\bar\lambda)\in\R^n\times\R^m$ is a local minimizer 
			of \eqref{eq:explicit_problem} for each $\bar\lambda\in K(\bar z)$, 
			and if $K$ is inner semicompact at $\bar z$ w.r.t.\ $\dom K$,
			then $\bar z$ is a local minimizer of \eqref{eq:implicit_problem}.
	\end{enumerate}
\end{theorem}
\begin{proof}
	The proof of the first assertion is, again, a simple consequence of the definition
	of the domain of a set-valued mapping. Hence, we merely prove the second statement.
	\\
	Suppose that $\bar z$ is not a local minimizer of \eqref{eq:implicit_problem}.
	Then we find a sequence $\{z_k\}_{k\in\N}\subset M\cap\dom K$ such that
	$z_k\to\bar z$ and $f(z_k)<f(\bar z)$ for all $k\in\N$.
	By definition of the domain, there is a sequence $\{\lambda_k\}_{k\in\N}\subset\R^m$ 
	such that $\lambda_k\in K(z_k)$ for all $k\in\N$,
	and noting that $K$ is inner semicompact at $\bar z$ w.r.t.\ $\dom K$,
	$\{\lambda_k\}_{k\in\N}$ may be chosen such that it possesses an accumulation point $\bar\lambda\in\R^m$.
	Closedness of $\gph K$ yields $(\bar z,\bar\lambda)\in\gph K$,
	i.e., $\bar\lambda\in K(\bar z)$.
	Then, however, $(\bar z,\bar\lambda)$ cannot be a local minimizer of
	\eqref{eq:explicit_problem}.
\end{proof}

In the literature, one can find diverse examples which illustrate that, in the second
statement of \cref{thm:local_minimizers}, one cannot abstain from postulating 
inner semicompactness of $K$ w.r.t.\ its domain or requiring local optimality 
for \emph{all} realizations of the implicit variable in \eqref{eq:explicit_problem}, 
see e.g.\ \cite[Example~4.8]{BenkoMehlitz2021}.
As it is clear from \cite{BenkoMehlitz2021},
and as we will illustrate in \cref{sec:consequences_examples}, the inner semicompactness
assumption is, in many reasonable applications, not restrictive.
However, the need to guarantee local minimality in \eqref{eq:explicit_problem}
for all potential choices of the implicit variable is a nontrivial requirement.
In practice, if it fails, it is typically violated for just one of these realizations,
see e.g.\ \cref{ex:MPVC_vs_MPCC} below.
This is very problematic from a conceptual point of view.
Typically, one aims to solve \eqref{eq:implicit_problem} by tackling \eqref{eq:explicit_problem}.
In situations where \eqref{eq:explicit_problem} is not convex,
a solution algorithm typically delivers 
a stationary point $(\bar z,\bar\lambda)\in\R^n\times\R^m$ of \eqref{eq:explicit_problem} which,
in well-behaved situations, is also a local minimizer of \eqref{eq:explicit_problem}. 
If $K(\bar z)$ is not a singleton, according to \cref{thm:local_minimizers}, 
this does not necessarily mean that $\bar z$
is a local minimizer of \eqref{eq:implicit_problem}. 

\cref{thm:local_minimizers} motivates the subsequently stated definition.
\begin{definition}\label{def:loc_min}
	We will refer to a global (local) minimizer of \eqref{eq:implicit_problem} 
	as an \emph{implicit} global (local) minimizer.
	Whenever $\bar z\in\R^n$ is a feasible point of \eqref{eq:implicit_problem}
	such that, for each $\bar\lambda\in K(\bar z)$, $(\bar z,\bar\lambda)$ is a global (local)
	minimizer of \eqref{eq:explicit_problem},
	we will call $\bar z$ an \emph{explicit} global (local) minimizer 
	of \eqref{eq:implicit_problem}.
\end{definition}

	\begin{remark}\label{rem:explicit_global_minimizer}
		Let us note that, according to \cref{thm:global_minimizers}, 
		a feasible point $\bar z\in\R^p$ of \eqref{eq:implicit_problem} 
		is an explicit global minimizer of the latter problem
		if any only if there exists $\bar\lambda\in K(\bar z)$
		such that $(\bar z,\bar\lambda)$ is a global minimizer
		of \eqref{eq:explicit_problem}.
		Beware that, in general, a similar equivalence does
		not hold true for explicit local minimizers as already 
		mentioned in the discussion following \cref{thm:local_minimizers}.
	\end{remark}

As a corollary of \cref{thm:local_minimizers}, we obtain the following result.
\begin{corollary}\label{cor:local_minimizers}
	\phantom{firstline}
	\begin{enumerate}
		\item Each implicit global (local) minimizer of \eqref{eq:implicit_problem} 
			is an explicit global (local) minimizer of \eqref{eq:implicit_problem}.
		\item Each explicit global minimizer of \eqref{eq:implicit_problem}
			is an implicit global minimizer of \eqref{eq:implicit_problem}.
		\item Each explicit local minimizer of \eqref{eq:implicit_problem} where
			$K$ is inner semicompact w.r.t.\ its domain is an implicit local minimizer 
			of \eqref{eq:implicit_problem}.
	\end{enumerate}
\end{corollary}

\subsection{Relations of stationary points}\label{sec:relations_stationary_points}

In this subsection, we exemplary focus on three concepts of stationarity
associated with \eqref{eq:implicit_problem} and, thus, \eqref{eq:explicit_problem}.
The first concept is a primal one and based on tangents to the feasible set
while the second and third concept are dual ones and make use of regular and limiting normals,
respectively.
Furthermore, we distinguish between implicit and explicit stationarity,
depending on whether it involves tangents/normals to the domain or the graph of $K$.

	Throughout the section,
	we will assume that the subsequently stated assumption holds.
	\begin{assumption}\label{ass:isc_at_reference_point}
		Let $\bar z\in\R^n$ be a feasible point of \eqref{eq:implicit_problem}
		where $K$ is inner semicompact w.r.t.\ $\dom K$.
	\end{assumption}

	As we have already mentioned in \cref{sec:variational_analysis},
	\cref{ass:isc_at_reference_point} guarantees that, locally around $\bar z$,
	$\dom K$ is closed.
	First, this ensures closedness 
	of the feasible set of \eqref{eq:implicit_problem}, locally around $\bar z$,
	which seems to be a rather standard assumption in mathematical optimization.
	Second, it allows for the computation of tangents and normals
	to $\dom K$ at $\bar z$ (that have been defined for locally closed sets here) 
	which is essential in the stationarity conditions
	addressing \eqref{eq:implicit_problem}. 
	Third, \cref{ass:isc_at_reference_point} puts us in position
	to apply \cref{lem:dom_graph_calculus} which will be essential for our analysis.
	Recall that \cref{ass:isc_at_reference_point} is rather mild
	and holds in many practically relevant scenarios
	as will be illustrated in \cref{sec:consequences_examples}.

Let us start with primal concepts of stationarity.
As we will rely on Bouligand's tangent cone and the associated graphical derivative,
these notions will be referred to as B-stationarity-type conditions.

\begin{definition}\label{def:primal_stationarity}
	Let \cref{ass:isc_at_reference_point} hold.
	\begin{enumerate}
		\item We refer to $\bar z$ as \emph{abstractly B-stationary}
			whenever we have
			\begin{equation}\label{eq:abstract_B_stat}
				\forall w\in T_{M\cap\dom K}(\bar z)\colon\quad
				f'(\bar z)w\geq 0.
			\end{equation}
		\item We refer to $\bar z$ as \emph{implicitly B-stationary}
			whenever we have
			\[
				\forall w\in T_{M}(\bar z)\cap T_{\dom K}(\bar z)\colon\quad
				f'(\bar z)w\geq 0.
			\]
		\item We refer to $\bar z$ as \emph{explicitly B-stationary}
			whenever
			\begin{equation}\label{eq:explicit_B_stat}
				\forall w\in T_M(\bar z)\cap \dom D K(\bar z,\bar\lambda)\colon\quad
				f'(\bar z)w\geq 0
			\end{equation}
			holds for all $\bar\lambda\in K(\bar z)$.
		\item Given $\bar\lambda\in K(\bar z)$, we refer to $\bar z$ as 
			\emph{explicitly B-stationary w.r.t.\ $\bar\lambda$}
			whenever \eqref{eq:explicit_B_stat} holds.
	\end{enumerate}
\end{definition}

The essence behind the definition of implicit B-stationarity is the separate
computation of tangents to $M$ and $\dom K$, i.e., to decouple the blocks of constraints
from each other in order to obtain conditions in terms of initial problem data. 
The abstract B-stationarity condition \eqref{eq:abstract_B_stat},
which is a necessary optimality condition for \eqref{eq:implicit_problem}
without any additional assumption,
is much more difficult to evaluate in practice 
as the geometry of the intersection $M\cap\dom K$ might be much more difficult
than the individual geometry of $M$ and $\dom K$, respectively.

Note that we can rewrite the condition \eqref{eq:explicit_B_stat} 
characterizing explicit B-stationarity as
\[
	\forall (w,d)\in T_{M\times\R^m}(\bar z,\bar\lambda)\cap T_{\gph K}(\bar z,\bar\lambda)
	\colon\quad
	f'(\bar z)w\geq 0,
\]
i.e., \eqref{eq:explicit_B_stat} can be interpreted as a B-stationarity-type condition 
associated with \eqref{eq:explicit_problem}.

Exploiting \cref{lem:dom_graph_calculus,lem:intersection_rule} as well as \cref{def:primal_stationarity}, 
the following result is immediate.

\begin{theorem}\label{thm:B_stat}
	Let \cref{ass:isc_at_reference_point} hold.
	Then the following assertions hold.
	\begin{enumerate}
		\item If $\bar z$ is an implicitly B-stationary point,
			it is abstractly B-stationary and explicitly B-stationary.
		\item The point $\bar z$ is an explicitly B-stationary point 
			if and only if it is an explicitly B-stationary point w.r.t.\ each $\bar\lambda\in K(\bar z)$.
		\item If $\bar z$ is an explicitly B-stationary point
			and if $K$ is inner calm* in the fuzzy sense at $\bar z$ w.r.t.\ $\dom K$,
			then $\bar z$ is an implicitly B-stationary point.
	\end{enumerate}
\end{theorem}

Next, we explain our interest in the stationarity notions 
which have been introduced in \cref{def:primal_stationarity}.
Therefore, let us introduce set-valued mappings 
$\Upsilon_\textup{imp}\colon\R^n\tto\R^n\times\R^n$
and $\Upsilon_\textup{exp}\colon\R^n\times\R^m\tto\R^n\times\R^n\times\R^m$ by means of
\begin{equation}\label{eq:feasibility_mappings}
	\begin{aligned}
		&\forall z\in\R^n\colon\quad&
			\Upsilon_\textup{imp}(z)&:=(z,z)-M\times\dom K,&
		\\
		&\forall (z,\lambda)\in\R^n\times\R^m\colon\quad&
			\Upsilon_\textup{exp}(z,\lambda)&:=(z,z,\lambda)-M\times\gph K.&
	\end{aligned}
\end{equation}
\begin{remark}\label{rem:abstract_to_implicit_B_stat}
		Let \cref{ass:isc_at_reference_point} hold.
		Furthermore, let $\bar z$ be an abstractly B-stationary point 
		of \eqref{eq:implicit_problem}.
	Furthermore, assume that $\Upsilon_\textup{imp}$ is metrically subregular
	at $(\bar z,(0,0))$. Then, due to \cref{lem:intersection_rule},
	$\bar z$ is implicitly B-stationary.
\end{remark}

\begin{theorem}\label{thm:B_stat_necessary}
	Let \cref{ass:isc_at_reference_point} hold.
	Then the following assertions hold.
	\begin{enumerate}
		\item If $\bar z$ is an implicit local minimizer of \eqref{eq:implicit_problem},
			then it is abstractly B-stationary.
		\item\label{item:implicit_B_stat} 
			If $\bar z$ is an implicit local minimizer of \eqref{eq:implicit_problem},
			and if $\Upsilon_\textup{imp}$ is metrically subregular at $(\bar z,(0,0))$,
			then $\bar z$ is implicitly B-stationary.
		\item If, for some $\bar\lambda\in K(\bar z)$,
			$(\bar z,\bar\lambda)$ is a local minimizer of \eqref{eq:explicit_problem},
			and if $\Upsilon_\textup{exp}$ is metrically subregular at $((\bar z,\bar\lambda),(0,0,0))$,
			then $\bar z$ is explicitly B-stationary w.r.t.\ $\bar\lambda$.
		\item\label{item:explicit_B_stat} 
			If $\bar z$ is an explicit local minimizer of \eqref{eq:implicit_problem},
			and if $\Upsilon_\textup{exp}$ is metrically subregular at $((\bar z,\bar\lambda),(0,0,0))$
			for each $\bar\lambda\in K(\bar z)$,
			then $\bar z$ is explicitly B-stationary.
	\end{enumerate}
\end{theorem}
\begin{proof}
	The first statement is a consequence of \cite[Theorem~6.12]{RockafellarWets1998},
	and the second assertion follows from the first one 
	with the aid of \cref{rem:abstract_to_implicit_B_stat}.
	To prove the third statement, we note that local optimality of $(\bar z,\bar\lambda)$
	for \eqref{eq:explicit_problem} yields
	\[
		\forall (w,d)\in T_{(M\times\R^m)\cap\gph K}(\bar z,\bar\lambda)\colon\quad
		f'(\bar z)w\geq 0
	\]
	by means of \cite[Theorem~6.12]{RockafellarWets1998}.
	Furthermore, it is easy to see that metric subregularity of $\Upsilon_\textup{exp}$
	at $((\bar z,\bar\lambda),(0,0,0))$ is equivalent to the metric subregularity of
	\[
		(z,\lambda)\mapsto(z,\lambda,z,\lambda)-M\times\R^m\times\gph K
	\]
	at $((\bar z,\bar\lambda),(0,0,0,0))$, so that the claim is a consequence of
	\cref{lem:intersection_rule}.
	The final assertion follows from the third one, 
	keeping \cref{def:loc_min} in mind.
\end{proof}

	As drafted at the end of \cref{sec:preliminaries},
	the metric subregularity assumptions on $\Upsilon_{\textup{imp}}$
	and $\Upsilon_{\textup{exp}}$ from \cref{thm:B_stat_necessary}
	are implied by validity of certain NNAMCQ-type qualification conditions
	which can be specified whenever variational descriptions of $M$ and $K$ are
	available.	

Let us now turn over to dual concepts of stationarity associated with \eqref{eq:implicit_problem}.
To start, we note that $\bar z\in\R^n$ being an abstractly B-stationary point of \eqref{eq:implicit_problem}
is equivalent to
\begin{equation}\label{eq:abstract_B_stat_dual}
	-\nabla f(\bar z)\in \widehat N_{M\cap\dom K}(\bar z)
\end{equation}
by the polar relationship between the tangent and the regular normal cone.
In order to develop this inclusion into more tractable conditions in terms of initial problem data,
i.e., variational objects in terms of $M$ and $K$, one has to face the commonly known fact that, variationally,
the regular normal cone behaves comparatively awkward in situations where the underlying sets are not convex or
at least not regular. However, depending on the problem structure, it still might be reasonable to rely on
\begin{equation}\label{eq:implicit_S_stationarity}
	-\nabla f(\bar z) \in \widehat N_M(\bar z)+\widehat N_{\dom K}(\bar z)
\end{equation}
as a necessary optimality condition for \eqref{eq:implicit_problem},
and keeping \cref{lem:dom_graph_calculus} in mind motivates the following definition.
\begin{definition}\label{def:dual_S_stationarity}
	Let \cref{ass:isc_at_reference_point} hold.
	\begin{enumerate}
		\item We refer to $\bar z$ as \emph{implicitly S-stationary} whenever \eqref{eq:implicit_S_stationarity} holds.
		\item We refer to $\bar z$ as \emph{explicitly S-stationary} whenever
			\begin{equation}\label{eq:explicit_S_stat}
				-\nabla f(\bar z)\in\widehat N_M(\bar z)+\widehat D^*K(\bar z,\bar\lambda)(0)
			\end{equation}
			holds for all $\bar\lambda\in K(\bar z)$.
		\item Given $\bar\lambda\in K(\bar z)$, we refer to $\bar z$ as \emph{explicitly S-stationary w.r.t.\ $\bar\lambda$}
			whenever \eqref{eq:explicit_S_stat} holds.
	\end{enumerate}
\end{definition}

Let us mention that the regular normal cone is associated with the concept of so-called strong stationarity
in the literature, see e.g.\ \cite{Mehlitz2020b} for the case of disjunctive optimization, 
and that is why we refer to the concepts in \cref{def:dual_S_stationarity} as S-stationarity-type conditions. 
In this regard, one might refer to \eqref{eq:abstract_B_stat_dual} as 
\emph{abstract S-stationarity} of $\bar z$, but the latter is equivalent
to abstract B-stationarity of $\bar z$
which is why we did not mention it in \cref{def:dual_S_stationarity}.
Similar to the definitions of B-stationarity, we decouple the blocks of constraints
in the definition of implicit and explicit S-stationarity in order to obtain
reasonable conditions in terms of problem data.
In contrast, the condition \eqref{eq:abstract_B_stat_dual} is difficult to evaluate exactly
as it requires the computation of the regular normal cone to an intersection.
Let us also mention that condition \eqref{eq:explicit_S_stat} 
characterizing explicit S-stationarity 
can be rewritten as
\[
	(-\nabla f(\bar z),0)
	\in
	\widehat N_{M\times\R^m}(\bar z,\bar\lambda) + \widehat N_{\gph K}(\bar z,\bar\lambda),
\] 
and, thus, can be interpreted as an S-stationarity-type condition 
associated with \eqref{eq:explicit_problem}.

To start our analysis of the concepts of S-stationarity, 
we first interrelate the stationarity notions in \cref{def:dual_S_stationarity}.
Similarly to \cref{thm:B_stat}, the following theorem is an immediate consequence of
\cref{lem:dom_graph_calculus,lem:intersection_rule}.

\begin{theorem}\label{thm:S_stat}
	Let \cref{ass:isc_at_reference_point} hold.
	Then the following assertions hold.
	\begin{enumerate}
		\item If $\bar z$ is an implicitly S-stationary point,
			it is abstractly B-stationary and explicitly S-stationary.
		\item The point $\bar z$ is explicitly S-stationary if and only if
			it is explicitly S-stationary w.r.t.\ each $\bar\lambda\in K(\bar z)$.
		\item\label{item:from_explicit_to_implicit_S_stat}
			If $\bar z$ is an explicitly S-stationary point,
			if $K(\bar z)$ is a singleton, 
			and if $K$ is inner calm* in the fuzzy sense at $\bar z$ w.r.t.\ $\dom K$,
			then $\bar z$ is an implicitly S-stationary point.
	\end{enumerate}
\end{theorem}

Let us mention that assuming $K(\bar z)$ to be a singleton is essential 
in \cref{thm:S_stat}\,\ref{item:from_explicit_to_implicit_S_stat}, 
see \cref{ex:explicit_S_stat_without_implicit_S_stat} below.
This additional assumption can be dropped if $\bar z$ is an interior
point of $M$ which yields $\widehat N_M(\bar z)=\{0\}$,
but that essentially means that the constraints $z\in M$ are locally
redundant around $\bar z$. The latter is a rare situation
which we are not going to investigate in the following.

We now mention situations where the concepts from \cref{def:dual_S_stationarity}
serve as necessary optimality conditions for \eqref{eq:implicit_problem}.
Similarly to \cref{thm:B_stat_necessary}, these results follow
from \cref{lem:intersection_rule} and \cref{def:loc_min}.

\begin{theorem}\label{thm:S_stat_necessary}
	Let \cref{ass:isc_at_reference_point} hold.
	Then the following assertions hold.
	\begin{enumerate}
		\item\label{item:implicit_S_stat} 
			If $\bar z$ is an implicit local minimizer of \eqref{eq:implicit_problem},
			and if there exists $\kappa>0$ such that
			\begin{equation}\label{eq:metric_inclusion_for implicit_S_stat}
				\widehat{\partial}\dist(\cdot,M\cap\dom K)(\bar z)
				\subset
				\kappa\bigl(
					\widehat\partial\dist(\cdot,M)(\bar z)
					+
					\widehat\partial\dist(\cdot,\dom K)(\bar z)
				\bigr),
			\end{equation}
			then $\bar z$ is implicitly S-stationary.
		\item If, for some $\bar\lambda\in K(\bar z)$,
			$(\bar z,\bar\lambda)$ is a local minimizer of \eqref{eq:explicit_problem},
			and if there exists $\kappa>0$ such that
			\begin{equation}\label{eq:metric_inclusion_for_explicit_S_stat}
				\begin{aligned}
				\widehat{\partial}\dist(\cdot,(M\times\R^m)&\cap\gph K)(\bar z,\bar\lambda)
				\\
				&\subset
				\kappa\bigl(
					\widehat\partial\dist(\cdot,M)(\bar z)\times\{0\}
					+
					\widehat\partial\dist(\cdot,\gph K)(\bar z,\bar\lambda)
				\bigr),
				\end{aligned}
			\end{equation}
			then $\bar z$ is explicitly S-stationary w.r.t.\ $\bar\lambda$.
		\item\label{item:explicit_S_stat} 
			If $\bar z$ is an explicit local minimizer of \eqref{eq:implicit_problem},
			and if, for each $\bar\lambda\in K(\bar z)$, 
			there exists $\kappa>0$ such that \eqref{eq:metric_inclusion_for_explicit_S_stat}
			is valid, then $\bar z$ is explicitly S-stationary.
	\end{enumerate}
\end{theorem}

	Following the arguments at the end of \cref{sec:preliminaries},
	inclusions \eqref{eq:metric_inclusion_for implicit_S_stat} 
	and \eqref{eq:metric_inclusion_for_explicit_S_stat} are implied
	by conditions of LICQ type.
	The latter ones can be worked out for many practically relevant realizations of $M$ and $K$. 

Recall from \cref{thm:B_stat_necessary} that abstract B-stationarity 
provides a necessary optimality condition for \eqref{eq:implicit_problem},
and that abstract B-stationarity is equivalent to \eqref{eq:abstract_B_stat_dual}.
To circumvent the limited available calculus of the regular normal cone, it is a standard approach 
to replace it by the larger limiting normal cone in \eqref{eq:abstract_B_stat_dual},
and this leads us to the subsequently stated definitions 
which are motivated by the calculus rules in \cref{lem:dom_graph_calculus,lem:intersection_rule}.

\begin{definition}\label{def:dual_stationarity}
	Let \cref{ass:isc_at_reference_point} hold.
	\begin{enumerate}
		\item We refer to $\bar z$ as \emph{abstractly M-stationary}
			whenever we have
			\[
				-\nabla f(\bar z) \in N_{M\cap\dom K}(\bar z).
			\]
		\item We refer to $\bar z$ as \emph{implicitly M-stationary}
			whenever we have
			\[
				-\nabla f(\bar z) \in N_M(\bar z) + N_{\dom K}(\bar z).
			\]
		\item We refer to $\bar z$ as \emph{explicitly M-stationary}
			whenever 
			\begin{equation}\label{eq:explicit_M_stat}
				- \nabla f(\bar z) \in N_M(\bar z) + D^*K(\bar z,\bar\lambda)(0)
			\end{equation}
			holds for all $\bar\lambda\in K(\bar z)$.
		\item Given $\bar\lambda\in K(\bar z)$, we refer to $\bar z$ as 
			\emph{explicitly M-stationary w.r.t.\ $\bar\lambda$}
			whenever \eqref{eq:explicit_M_stat} holds.
	\end{enumerate}
\end{definition}

As before, we note that explicit M-stationarity corresponds to an M-stationarity 
condition related to \eqref{eq:explicit_problem}.
Similar to \cref{thm:B_stat,thm:B_stat_necessary}, the following results are consequences
of \cref{lem:dom_graph_calculus,lem:intersection_rule}, 
\cref{def:loc_min,def:dual_stationarity},
and \cite[Theorem~6.12]{RockafellarWets1998}.
\begin{theorem}\label{thm:M_stat}
	Let \cref{ass:isc_at_reference_point} hold.
	Then the following assertions hold.
	\begin{enumerate}
		\item If $\bar z$ is abstractly M-stationary and if $\Upsilon_\textup{imp}$
			is metrically subregular at $(\bar z,(0,0))$,
			then $\bar z$ is implicitly M-stationary.
		\item If $\bar z$ is implicitly M-stationary,
			then there exists $\bar\lambda\in K(\bar z)$ such that
			$\bar z$ is explicitly M-stationary w.r.t.\ $\bar\lambda$.
		\item The point $\bar z$ is explicitly M-stationary if and only if
			it is explicitly M-stationary w.r.t.\ each $\bar\lambda\in K(\bar z)$.
	\end{enumerate}
\end{theorem}

\begin{theorem}\label{thm:M_stat_necessary}
	Let \cref{ass:isc_at_reference_point} hold.
	Then the following assertions hold.
	\begin{enumerate}
		\item If $\bar z$ is an implicit local minimizer of \eqref{eq:implicit_problem},
			then it is abstractly M-stationary.
		\item\label{item:implicit_M_stat} 
			If $\bar z$ is an implicit local minimizer of \eqref{eq:implicit_problem},
			and if $\Upsilon_\textup{imp}$ is metrically subregular at $(\bar z,(0,0))$,
			then $\bar z$ is implicitly M-stationary.
		\item If, for some $\bar\lambda\in K(\bar z)$,
			$(\bar z,\bar\lambda)$ is a local minimizer of \eqref{eq:explicit_problem},
			and if $\Upsilon_\textup{exp}$ is metrically subregular at $((\bar z,\bar\lambda),(0,0,0))$,
			then $\bar z$ is explicitly M-stationary w.r.t.\ $\bar\lambda$.
		\item\label{item:explicit_M-stat} 
			If $\bar z$ is an explicit local minimizer of \eqref{eq:implicit_problem},
			and if $\Upsilon_\textup{exp}$ is metrically subregular at $((\bar z,\bar\lambda),(0,0,0))$
			for each $\bar\lambda\in K(\bar z)$, then $\bar z$ is explicitly M-stationary.
	\end{enumerate}
\end{theorem}

It remains to interrelate the B-, S-, and M-stationarity concepts 
from \cref{def:primal_stationarity,def:dual_S_stationarity,def:dual_stationarity}
with each other.
\begin{theorem}\label{thm:relations_primal_dual_stat}
	Let \cref{ass:isc_at_reference_point} hold.
	Then the following assertions hold.
	\begin{enumerate}
		\item If $\bar z$ is an abstractly B-stationary point of \eqref{eq:implicit_problem},
			then it is abstractly M-stationary.
		\item If $\bar z$ is an implicitly S-stationary point of \eqref{eq:implicit_problem},
			then it is implicitly B- and M-stationary.
		\item If, for some $\bar\lambda\in K(\bar z)$, 
			$\bar z$ is an explicitly S-stationary point of \eqref{eq:implicit_S_stationarity}
			w.r.t.\ $\bar\lambda$,
			then it is explicitly B- and M-stationary w.r.t.\ $\bar\lambda$.
		\item If $\bar z$ is an explicitly S-stationary point of \eqref{eq:implicit_problem},
			then it is explicitly B- and M-stationary.
	\end{enumerate}
\end{theorem}
\begin{proof}
	In order to prove the first assertion,
	it is enough to recall that abstract B-stationarity is the same as \eqref{eq:abstract_B_stat_dual}
	as well as the fact that the regular normal cone to a set at some point is always a subset of the
	corresponding limiting normal cone.
	
	We proceed with the proof of the second assertion.
	Whenever $\bar z$ is implicitly S-stationary, it is trivially implicitly M-stationary
	due to the aforementioned general inclusion between the regular and limiting normal cone.
	Furthermore, implicit S-stationarity yields the existence of $\xi\in\widehat N_M(\bar z)$
	and $\chi\in\widehat N_{\dom K}(\bar z)$ such that $-\nabla f(\bar z)=\xi+\chi$.
	Hence, for each $w\in T_M(\bar z)\cap T_{\dom K}(\bar z)$, we find
	$f'(\bar z)w=-\xi^\top w-\chi^\top w\geq 0$,
	and the latter means that $\bar z$ is implicitly B-stationary.
	
	The third statement can be proven in analogous fashion as the second one, and
	the final assertion is an immediate consequence of the third one.
\end{proof}

	In the subsequently stated remark,
	we comment on the terminology used in \cref{def:primal_stationarity,def:dual_S_stationarity,def:dual_stationarity}.
	\begin{remark}\label{rem:explicit_stationarity}
		Let us note that the terminology of implicit and explicit stationarity,
		coined in \cref{def:primal_stationarity,def:dual_S_stationarity,def:dual_stationarity},
		has been chosen to underline the relationship to the implicit and explicit
		optimization problems \eqref{eq:implicit_problem} and \eqref{eq:explicit_problem},
		respectively.
		Furthermore, this terminology is consistent with the one used in the
		preceding paper \cite{BenkoMehlitz2021} where the concept
		of implicit variables has been introduced, see \cref{sec:comparison} as well.
		However, it should be observed that the computation or, at least, estimation 
		of the generalized derivatives of $K$, 
		which appear in the definitions of explicit B-, S-, and M-stationarity,
		might be quite involved and, potentially, requires additional assumptions
		in practically relevant scenarios.
		Hence, these stationarity conditions are not necessarily explicit in the sense
		that they are ready-to-use.
		To be fair, the same can be said about the implicit stationarity conditions
		as we typically face situations where $\dom K$ is a set of variationally challenging
		structure. 
	\end{remark}

Let us underline that our principle interest is in the identification of minimizers of \eqref{eq:implicit_problem},
and that \eqref{eq:explicit_problem} has to be interpreted as a seemingly more tractable surrogate problem.
Tackling \eqref{eq:implicit_problem} directly requires, exemplary, 
a subregularity assumption on $\Upsilon_\textup{imp}$
and yields implicit B- or M-stationarity of implicit local minimizers,
see \cref{thm:B_stat_necessary}\,\ref{item:implicit_B_stat} and \cref{thm:M_stat_necessary}\,\ref{item:implicit_M_stat}. 
Furthermore, following \cref{thm:local_minimizers} and our analysis in this subsection, 
meaningful necessary optimality conditions for \eqref{eq:implicit_problem} via \eqref{eq:explicit_problem}
should postulate explicit B- or M-stationarity
as these are necessary optimality conditions being guaranteed whenever $(\bar z,\bar\lambda)$ is a local
minimizer of \eqref{eq:explicit_problem} for each $\bar\lambda\in K(\bar z)$,
i.e., an explicit local minimizer of \eqref{eq:implicit_problem},
under suitable subregularity assumptions on $\Upsilon_\textup{exp}$,
see \cref{thm:B_stat_necessary}\,\ref{item:explicit_B_stat} and \cref{thm:M_stat_necessary}\,\ref{item:explicit_M-stat}.
Hence, one may ask about the strength of the associated constraint qualifications.
On the one hand, this would be metric subregularity of $\Upsilon_\textup{imp}$ at $(\bar z,(0,0))$
related to a direct treatment of \eqref{eq:implicit_problem}, and, on the other hand,
metric subregularity of $\Upsilon_\textup{exp}$ at $((\bar z,\bar\lambda),(0,0,0))$ for each $\bar\lambda\in K(\bar z)$,
corresponding to a treatment of \eqref{eq:explicit_problem}.
Subsequently, we compare the strength of these conditions.

\begin{proposition}\label{prop:comparison_subregularities}
	Let \cref{ass:isc_at_reference_point} hold.
	Furthermore, assume that, for each $\bar\lambda\in K(\bar z)$,
	$\Upsilon_\textup{exp}$ is metrically subregular at $((\bar z,\bar\lambda),(0,0,0))$.
	Then $\Upsilon_\textup{imp}$ is metrically subregular at $(\bar z,(0,0))$.
\end{proposition}
\begin{proof}
	We argue by contradiction.
	Suppose that $\Upsilon_\textup{imp}$ is not metrically subregular at $(\bar z,(0,0))$.
	Then, for each $k\in\N$, we find some $z_k\in\R^n$ such that
	\begin{align*}
			\norm{z_k-\bar z}
			&\geq
			\dist(z_k,M\cap\dom K)
			=
			\dist(z_k,\Upsilon_\textup{imp}^{-1}(0,0))
			\\
			&>
			k\,\dist((0,0),\Upsilon_\textup{imp}(z_k))
			=
			k\,(\dist(z_k,M)+\dist(z_k,\dom K))
	\end{align*}
	and $z_k\to \bar z$. Noting that $\dom K$ is closed locally around $\bar z$
	by the assumed inner semicompactness of $K$,
	for each sufficiently large $k\in\N$, we find $z_k'\in \dom K$ such that
	$\dist(z_k,\dom K)=\nnorm{z_k-z_k'}$.
	Furthermore, for any such $k\in\N$, the above yields
	\begin{equation}\label{eq:Upsilon_imp_not_ms}
		\norm{z_k-\bar z} 
		\geq
		\dist(z_k,M\cap \dom K)
		>
		k\,(\dist(z_k,M)+\nnorm{z_k-z_k'}).
	\end{equation}
	Particularly, $\norm{z_k-z_k'}\to 0$, i.e., $z_k'\to\bar z$.
	Exploiting the assumed inner semicompactness of $K$ again,
	there are a sequence $\{\lambda_k\}_{k\in\N}\subset\R^m$
	and a vector $\bar\lambda\in\R^m$
	such that $\lambda_k\in K(z_k')$ and $\lambda_k\to\bar\lambda$
	hold along a subsequence (without relabeling).
	The closedness of $\gph K$ guarantees $\bar\lambda\in K(\bar z)$.
	Due to $(z_k',\lambda_k)\in\gph K$, we have
	\[
		\dist((z_k,\lambda_k),\gph K)\leq \nnorm{z_k-z_k'}
	\]
	for each $k\in\N$, so that \eqref{eq:Upsilon_imp_not_ms} yields
	\begin{align*}
		\dist((z_k,\lambda_k),\Upsilon_\textup{exp}^{-1}(0,0,0))
		&=
		\dist((z_k,\lambda_k),(M\times\R^m)\cap\gph K)
		\\
		&\geq
		\dist(z_k,M\cap\dom K)
		\\
		&>
		k\,(\dist(z_k,M)+\nnorm{z_k-z_k'})
		\\
		&\geq
		k\,(\dist(z_k,M)+\dist((z_k,\lambda_k),\gph K))
		\\
		&=
		k\,\dist((0,0,0),\Upsilon_\textup{exp}(z_k,\lambda_k)).
	\end{align*}
	Consequently, $\Upsilon_\textup{exp}$ is not metrically subregular 
	at $((\bar z,\bar\lambda),(0,0,0))$ contradicting the postulated assumptions.
\end{proof}

The following example shows that the converse implication 
in \cref{prop:comparison_subregularities} does not hold. 
Hence, taking our above arguments into account,
\cref{prop:comparison_subregularities} shows that 
tackling \eqref{eq:implicit_problem} via \eqref{eq:explicit_problem} comes for the price
of more restrictive subregularity assumptions serving as constraint qualifications.

\begin{example}\label{ex:relation_of_subregularities}
	For $n:=m:=1$, let us consider $M:=\R_-$ and $K\colon\R\tto\R$ given by
	\[
		\forall z\in\R\colon\quad
		K(z)
		:=
		\begin{cases}
			\emptyset	&	\text{if }z<0,\\
			[-z^{1/2},z^{1/2}]	&	\text{if }z\geq 0.
		\end{cases}
	\]
	Obviously, we have $\dom K=\R_+$. 
	Hence, the associated mapping $\Upsilon_\textup{imp}$ is polyhedral
	and, thus, metrically subregular at each point of its graph.
	Let us now focus on the point $(\bar z,\bar\lambda):=(0,0)$
	which is the uniquely determined point in $(M\times\R)\cap\gph K$.
	One can easily check that
	\begin{align*}
		T_{(M\times\R)\cap\gph K}(\bar z,\bar\lambda)
		&=
		\{(0,0)\},
		\\
		T_{M\times\R}(\bar z,\bar\lambda)
		&=
		\R_-\times \R,
		\\
		T_{\gph K}(\bar z,\bar\lambda)
		&=
		\R_+\times\R,
		\\
		T_{M\times\R}(\bar z,\bar\lambda)\cap T_{\gph K}(\bar z,\bar\lambda)
		&=
		\{0\}\times\R
	\end{align*}
	which shows that the tangent cone intersection rule fails to hold with equality.
	Hence, the mapping $(z,\lambda)\mapsto(z,\lambda,z,\lambda)-M\times\R\times\gph K$ is not
	metrically subregular at $((\bar z,\bar\lambda),(0,0,0,0))$ due to
	\cref{lem:intersection_rule}, and this is equivalent to
	$\Upsilon_\textup{exp}$ failing to be metrically subregular at the point
	at $((\bar z,\bar\lambda),(0,0,0))$.
\end{example}

Above, metric subregularity of $\Upsilon_\textup{imp}$ and $\Upsilon_\textup{exp}$
may be seen as \emph{comparable} prototypical constraint qualifications for 
\eqref{eq:implicit_problem} and \eqref{eq:explicit_problem}, respectively,
and similar observations are likely to be made for other comparable concepts
of constraint qualifications which apply to \eqref{eq:implicit_problem} 
and \eqref{eq:explicit_problem} simultaneously.
Indeed, let us have a look at the qualification conditions used in \cref{thm:S_stat_necessary}
to guarantee implicit and explicit S-stationarity of a given local minimizer of \eqref{eq:explicit_problem}.

\begin{proposition}\label{prop:comparison_metric_inclusions}
	Let \cref{ass:isc_at_reference_point} hold.
	Furthermore, let $K$ be inner calm* in the fuzzy sense w.r.t.\ $\dom K$ at $\bar z$,
	and let $K(\bar z)$ be a singleton.
	Additionally, assume that, for the uniquely determined $\bar\lambda\in K(\bar z)$,
	condition \eqref{eq:metric_inclusion_for_explicit_S_stat} holds for some $\kappa>0$.
	Then \eqref{eq:metric_inclusion_for implicit_S_stat} holds for some (but not necessarily the same) $\kappa>0$.
\end{proposition}
\begin{proof}
	Due to \cref{lem:intersection_rule}, validity of \eqref{eq:metric_inclusion_for implicit_S_stat}
	for some $\kappa>0$ is equivalent to
	\begin{equation}\label{eq:regular_intersection_rule_implicit}
		\widehat N_{M\cap\dom K}(\bar z)
		=
		\widehat N_M(\bar z)+\widehat N_{\dom K}(\bar z),
	\end{equation}
	while validity of \eqref{eq:metric_inclusion_for_explicit_S_stat} for some $\kappa>0$ is equivalent to
	\begin{equation}\label{eq:regular_intersection_rule_explicit}
		\widehat N_{(M\times\R^m)\cap\gph K}(\bar z,\bar\lambda)
		=
		\widehat N_M(\bar z)\times\{0\}+\widehat N_{\gph K}(\bar z,\bar\lambda).
	\end{equation}
	We, thus, only need to show that the latter of these conditions implies the former.
	Therefore, we introduce a set-valued mapping $K_M\colon\R^n\tto\R^m$ by means of
	\[
		\forall z\in\R^n\colon\quad
		K_M(z)
		:=
		\begin{cases}
			K(z)	&	z\in M,\\
			\emptyset	&	z\notin M.
		\end{cases}
	\]
	As $M$ and $\gph K$ are closed, the same holds true for $\gph K_M=(M\times\R^m)\cap\gph K$.
	Furthermore, inner semicompactness of $K$ w.r.t. $\dom K$ at $\bar z$ implies
	inner semicompactness of $K_M$ w.r.t.\ $\dom K_M=M\cap\dom K$ at $\bar z$.
	Additionally, we have $K_M(\bar z)=\{\bar\lambda\}=K(\bar z)$.
	\Cref{lem:dom_graph_calculus}, thus, yields
	\begin{align*}
		\widehat N_{M\cap\dom K}(\bar z)
		=
		\widehat N_{\dom K_M}(\bar z)
		\subset
		\widehat D^* K_M(\bar z,\bar\lambda)(0).
	\end{align*}
	Hence, for $\zeta\in\widehat N_{M\cap \dom K}(\bar z)$,
	\eqref{eq:regular_intersection_rule_explicit} guarantees
	\begin{align*}
		(\zeta,0)
		&\in 
		\widehat N_{\gph K_M}(\bar z,\bar\lambda)
		=
		\widehat N_{(M\times\R^m)\cap\gph G}(\bar z,\bar\lambda)
		=
		\widehat N_M(\bar z)\times\{0\}+\widehat N_{\gph K}(\bar z,\bar\lambda).
	\end{align*}
	In other words, there is some $\xi\in\widehat N_M(\bar z)$ such that
	$\zeta-\xi\in\widehat D^* K(\bar z,\bar\lambda)(0)$.
	Now, inner calmness* in the fuzzy sense of $K$ w.r.t.\ $\dom K$ and \cref{lem:dom_graph_calculus}
	can be used to find $\zeta-\xi\in\widehat N_{\dom K}(\bar z)$,
	and this yields the inclusion $\subset$ in \eqref{eq:regular_intersection_rule_implicit}.
	Due to \cref{lem:intersection_rule}, the converse inclusion is generally valid.
	This completes the proof.
\end{proof}

Let us close with two additional comments.
First, one should note that validity of \eqref{eq:metric_inclusion_for implicit_S_stat} for some $\kappa>0$
(or, equivalently, \eqref{eq:regular_intersection_rule_implicit}) does not necessarily imply validity
of \eqref{eq:metric_inclusion_for_explicit_S_stat} for some (potentially different) $\kappa>0$
(or, equivalently, \eqref{eq:regular_intersection_rule_explicit}).
That can be easily seen in terms of \cref{ex:relation_of_subregularities}.
Second, as in \cref{thm:S_stat}, assuming that $K(\bar z)$ is a singleton is essential in \cref{prop:comparison_metric_inclusions},
see \cref{ex:explicit_S_stat_without_implicit_S_stat} below.

\subsection{A brief summery}\label{sec:summary}

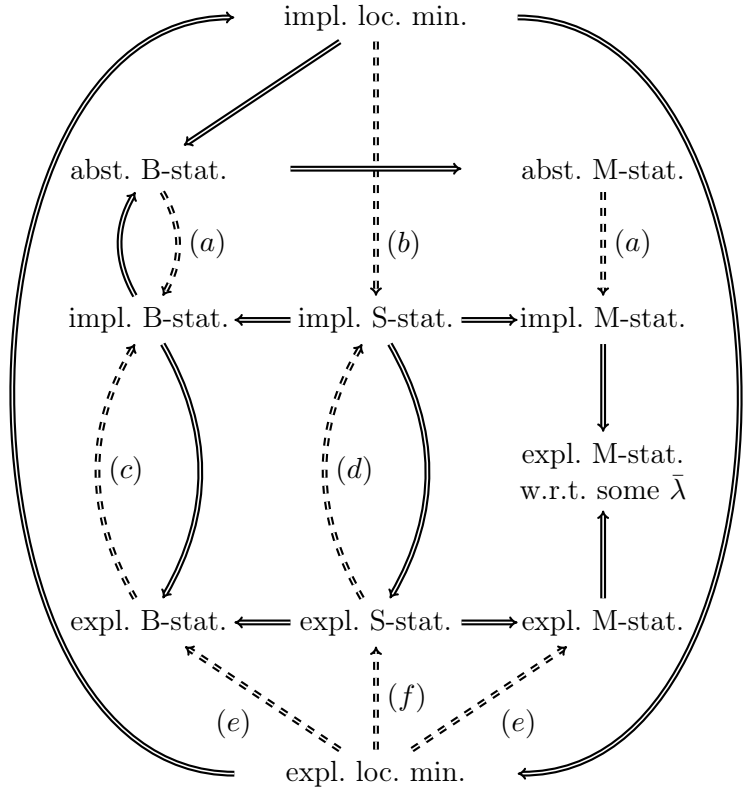
\begin{figure}[ht]
\centering
\begin{tikzpicture}[->]

  \node[punkt] at (0,0) 	(A){impl.\ loc.\ min.};
  \node[punkt] at (-3,-2) 	(B){abst.\ B-stat.};
  \node[punkt] at (3,-2) 	(C){abst.\ M-stat.};
  \node[punkt] at (-3,-4)   (D){impl.\ B-stat.};
  \node[punkt] at (0,-4) 	(E){impl.\ S-stat.};
  \node[punkt] at (3,-4) 	(F){impl.\ M-stat.};
  \node[punkt] at (3,-6)	(G){expl.\ M-stat.\ w.r.t.\ some $\bar\lambda$};
  \node[punkt] at (-3,-8)	(H){expl.\ B-stat.};
  \node[punkt] at (0,-8)	(I){expl.\ S-stat.};
  \node[punkt] at (3,-8)	(J){expl.\ M-stat.};
  \node[punkt] at (0,-10)	(K){expl.\ loc.\ min.};
  \node[right] at (0,-3) {$(b)$};

  \path     (A) edge[-implies,thick,double] node {}(B)
  			(A) edge[-implies,thick,double,dashed] node {}(E)
  			(A) edge[-implies,thick,double,bend angle=90,bend left] node {}(K)
            (B) edge[-implies,thick,double] node {}(C)
            (B) edge[-implies,thick,double,dashed,bend left] node[right] {$(a)$}(D)
            (C) edge[-implies,thick,double,dashed] node[right] {$(a)$}(F)
            (D) edge[-implies,thick,double,bend left] node {}(B)
            (D) edge[-implies,thick,double,bend left] node {}(H)
            (D) edge[-implies,thick,double] node {}(E)
            (F) edge[-implies,thick,double] node {}(E)
            (E) edge[-implies,thick,double,bend left] node {}(I)
			(F) edge[-implies,thick,double] node {}(G)
            (H) edge[-implies,thick,double,dashed,bend left] node[right] {$(c)$}(D)
            (H) edge[-implies,thick,double] node {}(I)
            (I) edge[-implies,thick,double,dashed,bend left] node[right] {$(d)$}(E)
            (J) edge[-implies,thick,double] node {}(G)
            (J) edge[-implies,thick,double] node {}(I)
            (K) edge[-implies,thick,double,bend angle=90,bend left] node {}(A)
            (K) edge[-implies,thick,double,dashed] node[below left] {$(e)$}(H)
            (K) edge[-implies,thick,double,dashed] node[right] {$(f)$}(I)
            (K) edge[-implies,thick,double,dashed] node[below right] {$(e)$}(J);
\end{tikzpicture}
\caption{
	Relations between properties of a feasible point $\bar z\in\R^n$ of \eqref{eq:implicit_problem}
	where $K$ is inner semicompact w.r.t.\ its domain.
	Dashed relations hold in the presence of additional assumptions outlined in
	\cref{sec:summary}.
}
\label{fig:summary}
\end{figure}

In \cref{fig:summary}, we summarize our findings from \cref{sec:relations_stationary_points}
in a clearly arranged graph which compares the \emph{implicit} and \emph{explicit} world.
Again we assume that \cref{ass:isc_at_reference_point} holds.
The mentioned additional assumptions used in \cref{fig:summary} are the following:
\begin{enumerate}
	\item $\Upsilon_\textup{imp}$ is metrically subregular at $(\bar z,(0,0))$;
	\item condition \eqref{eq:metric_inclusion_for implicit_S_stat} holds for some $\kappa>0$;
	\item $K$ is inner calm* in the fuzzy sense at $\bar z$ w.r.t.\ its domain;
	\item $K$ is inner calm* in the fuzzy sense at $\bar z$ w.r.t.\ its domain
			and $K(\bar z)$ is a singleton;
	\item for each $\bar\lambda\in K(\bar z)$,
			$\Upsilon_\textup{exp}$ is metrically subregular at $((\bar z,\bar\lambda),(0,0,0))$;
	\item for each $\bar\lambda\in K(\bar z)$,
			condition \eqref{eq:metric_inclusion_for_explicit_S_stat} holds for some $\kappa>0$.
\end{enumerate}

Below, we summarize the main messages which can be distilled from 
\cref{sec:minimizers,sec:relations_stationary_points} as well as \cref{fig:summary} so far.
\begin{itemize}
	\item Starting in the implicit setting, it is always possible to transfer any of the properties
		local minimality, B-stationarity, or S-stationarity
		into the explicit setting without further assumptions.
		When M-stationarity is under consideration, one has to be careful to note that
		an implicitly M-stationary point is, in general, only explicitly M-stationary
		w.r.t.\ some point in the associated image of $K$.
		In fact, it may happen that an implicitly M-stationary point is not explicitly M-stationary.
	\item Starting in the explicit setting, it is possible to transfer any of the properties
		local minimality, B-stationarity, or S-stationarity
		into the implicit setting for the price of additional assumptions.
		For local minimizers, one merely needs inner semicompactness of $K$ w.r.t.\ its domain which is a rather mild
		and reasonable assumption as it ensures the closedness of the feasible set of \eqref{eq:implicit_problem}.
		For B- and S-stationary, some quantitative version of inner semicompactness is needed
		for the transfer. 
		The situation is clearly different for explicit M-stationarity.
		We are not aware of any reasonable assumptions which guarantee that
		an explicitly M-stationary point is implicitly M-stationary
		even if the associated image of $K$ is a singleton.
	\item When solving \eqref{eq:implicit_problem} algorithmically, it is reasonable to expect an associated
		method to terminate at an implicitly M-stationary point,
		while algorithms treating the model \eqref{eq:explicit_problem} are likely to return points being
		explicitly M-stationary w.r.t.\ some point in the associated image of $K$.
		As outlined above, the latter might be weaker than the former with the gap being significant.
	\item The constraint qualifications needed to reasonably characterize an implicit local minimizer
		in terms of implicit B-, M-, or S-stationarity are less restrictive than the constraint qualifications
		one has to postulate in order to guarantee that an explicit local minimizer is explicitly
		B-, M-, or S-stationary, see \cref{prop:comparison_subregularities,prop:comparison_metric_inclusions}.
\end{itemize}
The above list indicates that treating implicit variables as explicit ones is likely to be disadvantageous.
The explicit model problem is likely to possess much more local minimizers and (in whatever sense)
stationary points than the implicit problem. Furthermore, the explicit model is likely to require more
restrictive constraint qualifications.
This strategy, which is frequently used to seemingly simplify complex structures within the feasibility conditions
of optimization problems, thus, comes at a very high price and should be avoided if possible.

\subsection{Comparison with a more complicated model}\label{sec:comparison}

In \cite{BenkoMehlitz2021}, the authors assume that the mapping $K$ is of the special type
\[
	\forall z\in\R^n\colon\quad
	K(z)
	:=
	\{\lambda\in F(z)\,|\,0\in G(z,\lambda)\},
\]
where $F\colon\R^n\tto\R^m$ and $G\colon\R^n\times\R^m\tto\R^s$ are given set-valued mappings
with a closed graph. 
This particular choice is motivated by several underlying applications
which can be modeled that way.
For further discussion, let us also introduce $H\colon\R^n\tto\R^s$ 
by means of
\[
	\forall z\in \R^n\colon\quad
	H(z):=\bigcup\limits_{\lambda\in F(z)}G(z,\lambda).
\]
The analysis in \cite{BenkoMehlitz2021} exploits $K$ as a so-called intermediate mapping
to handle the variational structure of $H$ which is a composition of set-valued mappings.
One can easily check that $H^{-1}(0)=\dom K$ holds true.
Hence, in this particular setting, 
\eqref{eq:implicit_problem} and \eqref{eq:explicit_problem} 
correspond to the implicit problem
\[
	\min\limits_z\{f(z)\,|\,z\in M,\,0\in H(z)\}
\]
and the associated explicit problem 
\[
	\min\limits_{z,\lambda}\{f(z)\,|\,z\in M,\,\lambda\in F(z),\,0\in G(z,\lambda)\}
\]
considered in \cite{BenkoMehlitz2021}.
Additionally, the concepts of implicit M-stationarity considered in \cref{def:dual_stationarity}
and \cite{BenkoMehlitz2021} coincide.
To evaluate the explicit concepts of M-stationarity from \cref{def:dual_stationarity},
one has to characterize the normal cone to $\gph K$ in this particular setting.
Noting that we have
\[
	\gph K
	=
	\{(z,\lambda)\in\R^n\times\R^m\,|\,((z,\lambda),((z,\lambda),0))\in\gph F\times\gph G\},
\]
it is possible to apply the pre-image rule from e.g.\ \cite[Section~5.1.4]{BenkoMehlitz2022}
in order to obtain, for $(\bar z,\bar\lambda)\in\gph K$,
\[
	N_{\gph K}(\bar z,\bar\lambda)
	\subset
	\left\{(\xi_1+\xi_2,\mu_1+\mu_2)\in\R^n\times\R^m\,\middle|\,
		\begin{aligned}
		&(\xi_1,\mu_1)\in N_{\gph F}(\bar z,\bar\lambda),
		\\
		&\exists \nu\in\R^s\colon\,			
			((\xi_2,\mu_2),-\nu)\in N_{\gph G}(\bar z,\bar\lambda)
		\end{aligned}
	\right\}
\]
whenever $(z,\lambda)\mapsto((z,\lambda),((z,\lambda),0))-\gph F\times\gph G$,
corresponding to the mapping 
$\mathfrak H\colon\R^n\times\R^m\tto\R^n\times\R^m\times\R^n\times\R^m\times\R^s$ 
defined in \cite[Section~4.2]{BenkoMehlitz2021},
is metrically subregular at $((\bar z,\bar\lambda),((0,0),((0,0),0)))$.
Using the definition of coderivatives, this would then yield the estimate
\[
	D^*K(\bar z,\bar\lambda)(0)
	\subset
	D^*F(\bar z,\bar\lambda)(\mu)
	+
	\{\xi\in\R^n\,|\,\exists\nu\in\R^s\colon\,(\xi,\mu)\in D^*G((\bar z,\bar\lambda),0)(\nu)\},
\]
showing that explicit M-stationarity of $\bar z$ w.r.t.\ $\bar\lambda\in K(\bar z)$
in the sense of \cref{def:dual_stationarity}
implies explicit M-stationarity 
of $\bar z$ w.r.t.\ $\bar\lambda$ used in \cite[Definition~4.10]{BenkoMehlitz2021}
in the presence of an additional subregularity assumption.
The latter is pretty much problem-tailored for this particular structure of $K$
while the former addresses the more general model problem \eqref{eq:implicit_problem}
which does not assume any underlying structure of $K$ apart from inner semicompactness
w.r.t.\ its domain.
Let us also note that metric subregularity of $\Upsilon_\textup{exp}$ at
$((\bar z,\bar\lambda),(0,0,0))$ and metric subregularity of
$\mathfrak H$ at $((\bar z,\bar\lambda),((0,0),((0,0),0)))$ together yield
metric subregularity of $(z,\lambda)\mapsto(z-M)\times \mathfrak H(z,\lambda)$
at $((\bar z,\bar\lambda),(0,((0,0),((0,0),0))))$
which is used as a qualification condition in \cite{BenkoMehlitz2021} to guarantee
explicit M-stationarity of a local minimizer $\bar z$ w.r.t.\ $\bar\lambda$ 
in the sense used therein.
To close, we want to point out that concepts similar to B- or S-stationarity 
from \cref{def:primal_stationarity,def:dual_S_stationarity}
are not discussed in \cite{BenkoMehlitz2021}.

\section{Selected consequences for the example problems}\label{sec:consequences_examples}

In this section,
we revisit the example problems from \cref{sec:examples} in the light of \cref{sec:analysis}.
For each of the example problems,
we show that the associated mapping $K$ used for the modeling
satisfies most of the technical assumptions which have been
exploited in \cref{sec:analysis}.
Furthermore, we illustrate some of the issues related to the appearance
of implicit variables by means of small academic examples. 

\subsection{Cardinality-constrained optimization}\label{sec:results_CCP}

In this subsection, we reinspect the cardinality-constrained optimization problem
\eqref{eq:CCP} and its surrogate model \eqref{eq:CCPref} in the light of our
findings from \cref{sec:analysis}.
Let us start by summarizing some important properties of the associated mapping $K_\textup{cc}$
from \eqref{eq:def_Kcc} in the subsequently stated lemma.
\begin{lemma}\label{lem:K_cc}
	\phantom{firstline}
	\begin{enumerate}
		\item The mapping $K_\textup{cc}$ possesses a closed graph.
		\item For $z\in\dom K_\textup{cc}$, $K_\textup{cc}(z)$ is a singleton
			if and only if $\norm{z}_0=\kappa$.
		\item The mapping $K_\textup{cc}$ possesses uniformly bounded image sets.
			Particularly, it is inner semicompact w.r.t.\ its domain at each point of its domain.
		\item The mapping $K_\textup{cc}$ is polyhedral.
			Particularly, it is metrically subregular at each point of its graph
			and inner calm* in the fuzzy sense w.r.t.\ its domain 
			at each point of its domain.
	\end{enumerate}
\end{lemma}
\begin{proof}
	The first statement is obvious as the set $\gph K_\textup{cc}$ is modeled
	with the aid of continuous functions.
	The second statement immediately follows by definition of $K_\textup{cc}$.
	The third assertion can be distilled from the fact that the images of $K_\textup{cc}$ 
	are included in the cube $[0,1]^n$.
	To verify the final statement, 
	we note that $\gph K_\textup{cc}$ possesses the representation
	\begin{equation}\label{eq:representation_gph_Kcc}
		\gph K_\textup{cc}
		=
		\bigcup\limits_{I\subset\{1,\ldots,n\}}
		\left\{(z,\lambda)\in\R^n\times\R^n\,\middle|\,
			\begin{aligned}
				&\mathtt e^\top\lambda\geq n-\kappa,\\
				&\forall i\in I\colon\, z_i=0,\,0\leq\lambda_i\leq 1,\\
				&\forall i\in\{1,\ldots,n\}\setminus I\colon\,\lambda_i=0
			\end{aligned}
		\right\}
	\end{equation}
	and, thus, is polyhedral.
	Hence, the final assertion follow from \cite[Theorem~3.4]{Benko2021}.
\end{proof}

\cref{lem:K_cc} shows that \eqref{eq:CCP} and the associated set-valued
mapping $K_\textup{cc}$ from \eqref{eq:def_Kcc} satisfy many of the
technical assumptions which have been exploited for the analysis
in \cref{sec:minimizers,sec:relations_stationary_points}.
Let us also mention that, whenever $M$ is the union of finitely many
convex polyhedral sets, then the associated mappings
$\Upsilon_\textup{imp}$ and $\Upsilon_\textup{exp}$ from
\eqref{eq:feasibility_mappings} are polyhedral and, thus,
metrically subregular at each point of their graphs
without further assumptions.

The upcoming example is included to illustrate that the assumption on
the set of implicit variables to be a singleton is essential 
in \cref{thm:S_stat} and \cref{prop:comparison_metric_inclusions}.

\begin{example}\label{ex:explicit_S_stat_without_implicit_S_stat}
	Let us consider the cardinality-constrained optimization problem
	\[
		\min\limits_z\{(z_1+1)^2+(z_2+1)^2\,|\,z_1-z_2=0,\,z\in\R^2_{\leq 1}\}
	\]
	at $\bar z:=0$ which is the uniquely determined feasible point of the problem and, thus,
	its global minimizer. 
	We find $K_\textup{cc}(\bar z)=\{\lambda\in\R^2\,|\,\lambda_1+\lambda_2\geq 1,\,0\leq \lambda \leq \mathtt e\}=:\Lambda$.
	Recall that $K_\textup{cc}$ is inner calm* in the fuzzy sense w.r.t.\ its domain at $\bar z$,
	see \cref{lem:K_cc}.
	\\
	First, we aim to show that $\bar z$ is explicitly but not implicitly S-stationary.
	Therefore, we have to compute $\widehat D^* K_\textup{cc}(\bar z,\bar\lambda)(0)$ 
	for each $\bar\lambda\in K_\textup{cc}(\bar z)$.
	Note that we can decompose $\gph K_\textup{cc}$ into three convex polyhedral sets as follows:
	\begin{equation*}
		\gph K_\textup{cc}
		=
		(\R\times\{(0,0,1)\})\cup(\{0\}\times\R\times\{(1,0)\})\cup(\{(0,0)\}\times\Lambda).
	\end{equation*}
	Hence, we can apply the union rule for the regular cone,
	see e.g.\ \cite[Lemma~2.2]{Mehlitz2020b}, to find
	\begin{equation}\label{eq:normal_cone_to_gph_Kcc}
		\widehat N_{\gph K_\textup{cc}}(\bar z,\bar\lambda)
		=
		\begin{cases}
			\R\times\{0\}\times\widehat N_\Lambda(\bar\lambda)  &	\bar\lambda=\mathtt e_1,
			\\
			\{0\}\times\R\times\widehat N_\Lambda(\bar\lambda)  &	\bar\lambda=\mathtt e_2,
			\\
			\R^2\times\widehat N_\Lambda(\bar \lambda)			&	\bar\lambda\in K_\textup{cc}(\bar x)\setminus\{\mathtt e_1,\mathtt e_2\},
		\end{cases}
	\end{equation}
	and this yields
	\[
		\widehat D^*K_\textup{cc}(\bar z,\bar\lambda)(0)
		=
		\begin{cases}
			\R\times\{0\}   &	\bar\lambda=\mathtt e_1,
			\\
			\{0\}\times\R   &	\bar\lambda=\mathtt e_2,
			\\
			\R^2			&	\bar\lambda\in K_\textup{cc}(\bar x)\setminus\{\mathtt e_1,\mathtt e_2\}.
		\end{cases} 
	\]
	Hence, for $\bar\lambda\in K_\textup{cc}(\bar z)\setminus\{\mathtt e_1,\mathtt e_2\}$,
	explicit S-stationarity w.r.t.\ $\bar\lambda$
	requires the existence of $\mu\in\R$ and $\zeta\in\R^2$ such that
	\begin{align*}
		\begin{pmatrix}
			0\\0
		\end{pmatrix}
		=
		\begin{pmatrix}
			2\\2
		\end{pmatrix}
		+
		\begin{pmatrix}
			1\\-1
		\end{pmatrix}
		\mu 
		+
		\zeta,	
	\end{align*}
	and this holds trivially.
	Furthermore, picking $\bar\lambda:=\mathtt e_i\in K_\textup{cc}(\bar z)$,
	explicit S-stationarity w.r.t.\ $\bar\lambda$ reduces to the existence of $\mu\in\R$ and $\zeta\in\R^2$ such that
	\begin{align*}
		\begin{pmatrix}
			0\\0
		\end{pmatrix}
		=
		\begin{pmatrix}
			2\\2
		\end{pmatrix}
		+
		\begin{pmatrix}
			1\\-1
		\end{pmatrix}
		\mu 
		+
		\zeta,
		\qquad
		\zeta_{3-i}=0.		
	\end{align*}
	For $i:=1$, one can choose $\mu:=2$, $\zeta_1:=-4$, and $\zeta_2:=0$ to solve this system,
	for $i:=2$, picking $\mu:=-2$, $\zeta_1:=0$, and $\zeta_2:=-4$ provides a solution.
	Summing this up, $\bar z$ is explicitly S-stationary.
	Due to $\widehat N_{\R^2_{\leq 1}}(\bar z)=\{0\}$, see e.g.\ \cite[Lemma~2.3]{PanXiuFan2017},
	implicit S-stationarity reduces to the existence of $\mu\in\R$ such that
	\begin{align*}
		\begin{pmatrix}
			0\\0
		\end{pmatrix}
		=
		\begin{pmatrix}
			2\\2
		\end{pmatrix}
		+
		\begin{pmatrix}
			1\\-1
		\end{pmatrix}
		\mu 
		,		
	\end{align*}
	but this condition does not hold. Hence, $\bar z$ is not implicitly S-stationary.
	\\
	Next, we aim to illustrate that the qualification condition \eqref{eq:regular_intersection_rule_explicit} holds for
	each $\bar\lambda\in K_\textup{cc}(\bar z)$ while \eqref{eq:regular_intersection_rule_implicit} is violated.
	Recalling \cref{lem:intersection_rule}, this illustrates that we cannot avoid to assume that the set
	of implicit variables is a singleton in \cref{prop:comparison_metric_inclusions}.
	Violation of \eqref{eq:regular_intersection_rule_implicit} is obvious 
	as we have shown above that $\bar z$ is not implicitly S-stationary despite being a
	minimizer of the problem of interest,
	see \cref{lem:intersection_rule} and \cref{thm:S_stat_necessary}.
	Hence, let us show validity of \eqref{eq:regular_intersection_rule_explicit} for each $\bar\lambda\in K_\textup{cc}(\bar z)$.
	Therefore, let us set $M:=\{z\in\R^2\,|\,z_1-z_2=0\}$.
	On the one hand, we find
	\begin{align*}
		(M\times\R^2)\cap\gph K_\textup{cc}
		&=
		\{(t,t,\lambda_1,\lambda_2)\in\R^4\,|\,\lambda_1+\lambda_2\geq 1,\,0\leq \lambda \leq \mathtt e,\,\lambda_1t=\lambda_2t=0\}
		\\
		&=
		\{(0,0)\}\times\Lambda,
	\end{align*}
	i.e.,
	\begin{align*}
		\widehat N_{(M\times\R^2)\cap\gph K_\textup{cc}}(\bar z,\bar\lambda)=\R^2\times\widehat N_{\Lambda}(\bar\lambda)
	\end{align*}
	for each $\bar\lambda\in K_\textup{cc}(\bar z)$.
	Picking $\bar\lambda\in K_\textup{cc}(\bar z)\setminus\{\mathtt e_1,\mathtt e_2\}$
	and relying on \eqref{eq:normal_cone_to_gph_Kcc},
	\eqref{eq:regular_intersection_rule_explicit} holds trivially as we have
	\begin{align*}
		\widehat N_{M\times\R^2}(\bar z,\bar\lambda)+\widehat N_{\gph K_\textup{cc}}(\bar z,\bar\lambda)
		&\supset
		\widehat N_{\gph K_\textup{cc}}(\bar z,\bar\lambda)
		\\
		&=
		\R^2\times\widehat N_{\Lambda}(\bar\lambda)
		=
		\widehat N_{(M\times\R^2)\cap\gph K_\textup{cc}}(\bar z,\bar\lambda),
	\end{align*}
	and the converse inclusion is generally valid, see \cref{lem:intersection_rule}.
	Choosing $\bar\lambda:=\mathtt e_1$, \eqref{eq:normal_cone_to_gph_Kcc} 
	and $\widehat N_M(\bar z)=\{\zeta\in\R^2\,|\,\zeta_1+\zeta_2=0\}$ yield
	\begin{align*}
		\widehat N_M(\bar z)\times\{(0,0)\}+\widehat N_{\gph K_\textup{cc}}(\bar z,\bar\lambda)
		=
		\R^2\times\widehat N_{\Lambda}(\bar\lambda)
		=
		\widehat N_{(M\times\R^2)\cap\gph K_\textup{cc}}(\bar z,\bar\lambda).
	\end{align*}
	Hence, \eqref{eq:regular_intersection_rule_explicit} holds for $\bar\lambda:=\mathtt e_1$,
	and in analogous fashion, one can verify validity of this condition for $\bar\lambda:=\mathtt e_2$.
	Thus, \eqref{eq:regular_intersection_rule_explicit} holds true for each $\bar\lambda\in K_\textup{cc}(\bar z)$.
\end{example}

In order to make the implicit notions of B-, S-, and M-stationarity associated with
\eqref{eq:CCP} practically useful, the tangent and normal cones to $\R^n_{\leq\kappa}$
have to be computed. Precise formulas are already available and can be found, e.g.,
in \cite[Lemma~2.3]{PanXiuFan2017}.
Let us also mention that the representation \eqref{eq:representation_gph_Kcc}
can be used to exactly compute the graphical derivative as well as the regular coderivative
of $K_\textup{cc}$ via the union rule from \cite[Lemma~2.2]{Mehlitz2020b} which already has been exploited 
in \cref{ex:explicit_S_stat_without_implicit_S_stat}.
Similarly, \cite[Corollary~1]{AdamCervinkaPistek2016} can be used to exactly compute the
limiting coderivative of $K_\textup{cc}$.
Hence, ready-to-use formulas for the explicit B-, S-, and M-stationarity notions are available.
In order to avoid prolonging the paper unnecessarily,
we leave these calculations to the interested reader.
Let us, however, note that the arising explicit stationarity notions of S- and M-type 
correspond to the ones which have been introduced for \eqref{eq:CCPref} 
in \cite[Definition~4.6]{BurdakovKanzowSchwartz2016}.

\subsection{Optimization problems with vanishing constraints}\label{sec:MFVC_consequences}

In this subsection, we readdress the cardinality-constrained optimization problem \eqref{eq:MPVC}
and its surrogate model \eqref{eq:MPVCref} 
while keeping our findings from \cref{sec:analysis} in mind.
Again, we start our investigations by summarizing essential properties 
of the associated mapping $K_\textup{vc}$ from \eqref{eq:def_Kvc}.

\begin{lemma}\label{lem:K_vc}
	\phantom{firstline}
	\begin{enumerate}
		\item The mapping $K_\textup{vc}$ possesses a closed graph.
		\item Fix $z\in\dom K_\textup{vc}$.
			Then $K_\textup{vc}(z)$ is a singleton
			if and only if $H_i(z)>0$ holds for all $i\in\{1,\ldots,m\}$.
			Otherwise, $K_\textup{vc}(z)$ is unbounded.
		\item The mapping $K_\textup{vc}$ is inner semicompact 
			w.r.t.\ its domain at each point of its domain.
		\item The mapping $K_\textup{vc}$ is inner calm* in the fuzzy sense 
			w.r.t.\ its domain at each point of its domain.
	\end{enumerate}
\end{lemma}
\begin{proof}
	The first assertion is trivial as $\gph K_\textup{vc}$ is modeled
	via continuous functions.
	The second statement is an immediate consequence of the representation
	\begin{equation}\label{eq:representation_Kvc}
		\forall z\in\dom K_\textup{vc}\colon\quad
		K_\textup{vc}(z)
		=
		\left\{
			\lambda\in\R^m\,\middle|\,
			\begin{aligned}
				&\lambda_i\geq \max(G_i(z),0)	&	&\forall i\in I^0_{H}(z)
				\\
				&\lambda_i=0				&	&\forall i\in I^+_H(z)
			\end{aligned}
		\right\},
	\end{equation}
	where we used
	\[	
		I^0_H(z):=\{i\in\{1,\ldots,m\}\,|\,H_i(z)=0\},
		\qquad
		I^+_H(z):=\{i\in\{1,\ldots,m\}\,|\,H_i(z)>0\}.
	\]
	In order to prove the final two assertions of the lemma,
	we introduce $\psi\colon\dom K_\textup{vc}\to\R^m$ by means of
	$\psi(z):=\max(G(z),0)$ for each $z\in\dom K_\textup{vc}$.
	Observe that $\psi(z)\in K_\textup{vc}(z)$ holds for each
	$z\in\dom K_\textup{vc}$.
	As $G$ is continuous, $\psi$ is continuous as well,
	and this shows the claimed inner semicompactness of $K_\textup{vc}$
	w.r.t.\ its domain at each point of its domain.
	Furthermore, as $G$ is even continuously differentiable,
	it is locally Lipschitz continuous,
	and $\psi$, thus, possesses the same property.
	This is already enough to guarantee that $K_\textup{vc}$
	is inner calm* in the fuzzy sense w.r.t.\ its domain at each
	point of its domain.
\end{proof}

The second assertion of \cref{lem:K_vc} underlines that
$K_\textup{vc}$ is, in nontrivial situations, not
locally bounded at the points of its domain,
but this does not effect the inherent inner semicompactness
w.r.t.\ its domain which has been proven.

\cref{lem:K_vc} illustrates that most of the technical assumptions
which have been exploited in \cref{sec:minimizers,sec:relations_stationary_points}
hold for the mapping $K_\textup{vc}$ from \eqref{eq:def_Kvc}.
Whenever the functions $G$ and $H$ are affine, $K_\textup{vc}$ is a polyhedral mapping.
If, additionally, $M$ is the union of finitely many convex polyhedral sets,
the associated mappings $\Upsilon_\textup{imp}$ and $\Upsilon_\textup{exp}$ 
from \eqref{eq:feasibility_mappings} are polyhedral and, thus, metrically subregular at each point
of their graphs without further assumptions.

In \cite[Lemma~1]{AchtzigerKanzow2008}, the authors claim that \eqref{eq:MPVC} and
\eqref{eq:MPVCref} are equivalent w.r.t.\ local minimizers in a certain sense. This
assertion is, unluckily, not true as the following simple example indicates.

\begin{example}\label{ex:MPVC_vs_MPCC}
	Let us consider the simple vanishing-constrained problem
	\begin{equation}\label{eq:simple_MPVC}
		\min\limits_z\{-z_2\,|\,z_2\geq 0,\,z_1z_2\leq 0\}
	\end{equation}
	as well as its complementarity-constrained reformulation
	\begin{equation}\label{eq:simple_MPCC}
		\min\limits_{z,\lambda}\{-z_2\,|\,z_1-\lambda\leq 0,\,z_2\geq 0,\,\lambda\geq 0,\,z_2\lambda=0\}.
	\end{equation}
	We note that $(\bar z,\bar \lambda):=((0,0),1)$ is a feasible point
	of \eqref{eq:simple_MPCC} which is a local minimizer.
	Indeed, small perturbations of this point in the feasible set always require $z_2=0$
	which explains this local minimality.
	However, $\bar z$ is not a local minimizer of \eqref{eq:simple_MPVC}.
	This refutes \cite[Lemma~1]{AchtzigerKanzow2008}.
	
	Let us explain what goes wrong here.
	Therefore, we make use of \cref{thm:local_minimizers}.
	Exploiting \eqref{eq:representation_Kvc}, $K_\textup{vc}(\bar z)=\R_+$ is valid,
	and it is easy to check that
	$(\bar z,\tilde\lambda)$ for $\tilde\lambda:=0$ is not a local
	minimizer of \eqref{eq:simple_MPCC}.
	Hence, \cref{thm:local_minimizers} indeed shows that $\bar z$ cannot be
	a local minimizer of \eqref{eq:simple_MPVC}.
\end{example}

In order to make the implicit and explicit notions of B-, S-, and M-stationarity
for \eqref{eq:MPVC} practically accessible, tangent and normal cones to the domain and graph
of $K_\textup{cc}$ have to be computed.
This can be done with the aid of suitable pre-images rules,
see e.g.\ \cite[Theorems~6.14 and~6.31]{RockafellarWets1998}.
For the implicit notions, 
we end up with the concepts promoted in \cite{Hoheisel2009},
while the explicit notions lead to the well-known stationarity concepts
from complementarity-constrained optimization, see e.g.\ \cite{Ye2005} for an overview,
under mild qualification conditions.
	Following our arguments at the end of \cref{sec:preliminaries},
	metric subregularity of $\Upsilon_\textup{imp}$ is implied
	by the constraint qualification MPVC-NNAMCQ,
	see \cite[Definition~2]{MokhtavayiHeydariKanzi20201},
	while \eqref{eq:metric_inclusion_for implicit_S_stat} in 
	valid in the presence of so-called MPVC-LICQ,
	see \cite[Definition~5.1.1]{Hoheisel2009}.
	Furthermore, metric subregularity of $\Upsilon_\textup{exp}$
	or \eqref{eq:metric_inclusion_for_explicit_S_stat}
	holds whenever MPCC-NNAMCQ or MPCC-LICQ,
	see \cite[Definitions~2.8 and~2.10]{Ye2005},
	is valid for the complementarity-constrained surrogate problem
	\eqref{eq:MPVCref}, respectively.

In \cite{AchtzigerKanzow2008}, the authors already point out that the consideration
of the surrogate model \eqref{eq:MPVCref} instead of \eqref{eq:MPVC} might not be
recommendable due to the larger dimension of \eqref{eq:MPVCref} and the potential
set-valuedness of the mapping $K_\textup{vc}$ although no evidence is presented.
The results in \cref{sec:analysis} as well as \cite{BenkoMehlitz2021}, however, 
underline this conjecture on the theoretical side as well.

\subsection{Bilevel optimization}

To close the section,
we reinspect the bilevel optimization problem \eqref{eq:BOP}
and its associated reformulation \eqref{eq:BOPref}
in the light of \cref{sec:analysis}.
Similar to the previous subsections,
we would like to start by summarizing the most important properties 
of $K_\textup{bop}$ from \eqref{eq:def_Kbop}.
Recall our standing assumption that \eqref{eq:lower_level}
is a convex optimization problem that
satisfies Slater's constraint qualification for each $z_1\in\dom\Psi$.

\begin{lemma}\label{lem:K_bop}
	\phantom{firstline}
	\begin{enumerate}
		\item The mapping $K_\textup{bop}$ possesses a closed graph.
		\item Fix $z\in\dom K_\textup{bop}$.
			Then $K_\textup{bop}(z)$ is a singleton
			if and only if, for $\lambda\in K_\textup{bop}(z)$, 
			\begin{equation}\label{eq:SMFC}
				\left.
				\begin{aligned}
				&0=\sum\nolimits_{i\in I(z)}\mu_i\nabla_{z_2}G_i(z),
				\\
				&\forall i\in I_0(z,\lambda)\colon\,\mu_i\geq 0
				\end{aligned}
				\right\}
				\quad\Longrightarrow\quad
				\mu_i=0\quad \forall i\in I(z)
			\end{equation}
			where
			\[
				I(z):=\{i\in\{1,\ldots,m\}\,|\,G_i(z)=0\},
				\qquad
				I_0(z,\lambda):=\{i\in I(z)\,|\,\lambda_i=0\}.
			\]
		\item The mapping $K_\textup{bop}$ is locally bounded at each point of its domain.
			Particularly, $K_\textup{bop}$ is inner semicompact w.r.t.\ its domain
			at each point of its domain.
		\item Fix $z\in\dom K_\textup{bop}$.
			Then the mapping $K_\textup{bop}$ is inner calm* in the fuzzy sense 
			w.r.t.\ its domain at $z$ if one of the following conditions holds:
			\begin{enumerate}[label=(\roman*)]
				\item $G$ does not depend on $z_1$,
				\item for $\lambda\in K_\textup{bop}(z)$,
						\eqref{eq:SMFC} is valid.
			\end{enumerate}
	\end{enumerate}
\end{lemma}
\begin{proof}
	The first assertion follows from continuous differentiability of $g$ and $G$.
	The second assertion can be found in \cite[Section~5]{Wachsmuth2013}.
	In the third assertion,
	uniform boundedness of the images of $K_\textup{bop}$ 
	locally around each point of its domain can be distilled
	from e.g.\ \cite[Proposition~4.43]{BonnansShapiro2000}, 
	and this implies that $K_\textup{bop}$ is inner semicompact w.r.t.\ its domain 
	at each point of its domain.
	The final assertion follows from \cite[Theorem~3.9]{Benko2021}
	and \cite[Proposition~4.47]{BonnansShapiro2000}.
\end{proof}

Let us mention that \eqref{eq:SMFC} is referred to as
strict Mangasarian--Fromovitz condition in the literature.
The latter is not a constraint qualification in the narrower sense
as it depends implicitly on the objective function $g$ via the multiplier $\lambda$.
Note that it holds in the presence of the linear independence constraint qualification.

It should be noted that, in order to guarantee inner semicompactness of $K_\textup{bop}$
w.r.t.\ its domain, one does not necessarily need to assume 
validity of Slater's constraint qualification
which already yields uniform boundedness of the images of $K_\textup{bop}$
locally around each point of its domain.
In \cite[Proposition~3.2, Theorem~3.3]{GfrererMordukhovich2017},
it is shown that so-called Robinson stability
(which is, in general, weaker than Slater's constraint qualification)
or a parametric version of the constant rank constraint qualification
(which is independent of Slater's constraint qualification)
provide sufficient conditions for the inner semicompactness of $K_\textup{bop}$
w.r.t.\ its domain. 
Observe that, for reformulation \eqref{eq:BOPref} of \eqref{eq:BOP} to be reasonable,
one merely requires validity of some constraint qualification 
for \eqref{eq:lower_level} for all $z_2\in\Psi(z_1)$ and $z_1\in\dom\Psi$,
i.e., one can indeed relax the inherent validity of Slater's constraint qualification 
if necessary.

\cref{lem:K_bop} illustrates that most of the technical assumptions which have been exploited
for the analysis in \cref{sec:minimizers,sec:relations_stationary_points} 
can be guaranteed for the mapping $K_\textup{bop}$ from \eqref{eq:def_Kbop}. 
Whenever $M$ is the union of finitely many convex polyhedral sets and \eqref{eq:lower_level}
takes the form
\begin{equation}\label{eq:quadratic_lower_level}
	\min\limits_{z_2}\{\tfrac12 z_2^\top Qz_2+(Cz_1)^\top z_2\,|\, Az_1+Bz_2\leq b\}
\end{equation}
for a positive semidefinite matrix $Q\in\R^{n_2\times n_2}$,
matrices $A\in\R^{m\times n_1}$, $B\in\R^{m\times n_2}$, and $C\in\R^{n_2\times n_1}$,
as well as a vector $b\in\R^m$,
then $\Upsilon_\textup{imp}$ and $\Upsilon_\textup{exp}$ 
from \eqref{eq:feasibility_mappings} are polyhedral mappings and, thus,
metrically subregular at each point of their graphs without further assumptions.
Indeed, in this case, $K_\textup{bop}$ is a polyhedral mapping and, thus,
its graph and domain are the union of finitely many convex polyhedral sets.
Note that we do not need to assume validity of Slater's constraint qualification
for \eqref{eq:quadratic_lower_level} as its constraints are affine so that
Abadie's constraint qualification holds everywhere.

Applying the results from \cref{sec:minimizers} to \eqref{eq:BOP} and $K_\textup{bop}$
from \eqref{eq:def_Kbop} restores findings from \cite{DempeDutta2012}.
	More precisely,
	\cref{thm:global_minimizers} recovers \cite[Theorems~2.1 and~2.3]{DempeDutta2012}
	while \cref{thm:local_minimizers} covers \cite[Theorem~3.2]{DempeDutta2012}
	when keeping \cref{lem:K_bop} in mind.
Metric-subregularity-type constraint qualifications for \eqref{eq:BOP}
and its reformulation \eqref{eq:BOPref} have been compared in \cite{AdamHenrionOutrata2018}
and came up with the same conclusions as \cref{prop:comparison_subregularities}.
To apply the implicit concepts of B-, S-, and M-stationarity 
from \cref{sec:relations_stationary_points} to the setting at hand, 
a tractable reformulation of $\dom K_\textup{bop}=\gph\Psi$ has to be chosen.
Let us mention two options to do this.

First, exploiting the so-called optimal value function 
$\varphi\colon\R^{n_1}\to\R\cup\{-\infty,\infty\}$ of \eqref{eq:lower_level} given by
\[
	\forall z_1\in\R^{n_1}\colon\quad
	\varphi(z_1):=\inf\limits_{z_2}\{g(z)\,|\,G(z)\leq 0\},
\]
we find
\[
	\gph\Psi=\{z\in\R^{n}\,|\,g(z)-\varphi(z_1)\leq 0,\,G(z)\leq 0\}.
\]
Tangents and normals to this set can be computed with the aid of a suitable pre-image rule. 
One has to be aware of two difficulties which have to be faced when using this approach.
On the one hand, $\varphi$ is likely to be nonsmooth.
On the other hand, any nonsmooth version of the Mangasarian--Fromovitz constraint qualification
fails at all points of $\gph\Psi$ when modeled in this way,
see \cite[Proposition~3.2]{YeZhu1995}. 
Hence, comparatively weak qualification conditions have
to be employed in order to guarantee metric subregularity of $\Upsilon_\textup{imp}$
in the present setting,
e.g.\ the assumptions of \cite[Lemma~5.1]{JolaosoMehlitzZemkoho2025}.
A theoretical and numerical comparison of this approach with the one promoted via
\eqref{eq:BOPref} is provided in \cite{ZemkohoZhou2021} 
and comes up with the recommendation to rely on \eqref{eq:BOP} whenever possible.

Second, one could exploit the feasibility mapping $\Gamma\colon\R^{n_1}\tto\R^{n_2}$
associated with \eqref{eq:lower_level} and given by 
$\Gamma(z_1):=\{z_2\in\R^{n_2}\,|\,G(z)\leq 0\}$, $z_1\in\R^{n_1}$,
in order to state
\[
	\gph\Psi
	=
	\{z\in\R^{n}\,|\,-\nabla_{z_2}g(z)\in \widehat N_{\Gamma(z_1)}(z_2)\}.
\]
Then the computation of tangents and normals to $\dom K_\textup{bop}$ requires
the computation of tangents and normals to the graph of the normal cone mapping
$z\mapsto \widehat N_{\Gamma(z_1)}(z_2)$,
but that is a well-known topic in variational analysis nowadays,
see e.g.\ \cite{GfrererOutrata2016a,GfrererOutrata2016b} for the case where $G$
does not depend on $z_1$
and \cite{BenkoGfrererOutrata2020,GfrererOutrata2020} for the more general situation.
These findings require twice continuous differentiability of $g$ and $G$.
	Let us note that the considerations in \cite{AdamHenrionOutrata2018}
	are based on this representation of $\gph\Psi$.

Notions of explicit B-, S-, and M-stationarity type for \eqref{eq:BOP} 
turn out to equal B-, S-, and M-stationarity associated with the complementarity-constrained
problem \eqref{eq:BOPref} 
and provide necessary optimality conditions
in the presence of reasonable qualification conditions
like MPCC-NNAMCQ and MPCC-LICQ
which we already mention in \cref{sec:MFVC_consequences}.
Let us note that, due to the special structure of \eqref{eq:BOPref},
both constraint qualifications are likely to fail
if $K_\textup{bop}$ is not singleton-valued at the reference point,
see \cite[Section~3]{GfrererYe2017}.

\section{Implicit variables in numerical cardinality-constrained optimization}\label{sec:numerics}

In this section, we aim to illustrate issues related to implicit variables
by means of a computational experiment.
Here, we focus on cardinality-constrained optimization problems.
For a fair numerical comparison of problem
\eqref{eq:CCP} and its reformulation \eqref{eq:CCPref},
we will rely on the augmented Lagrangian framework from \cite{JiaKanzowMehlitzWachsmuth2023} 
that nicely applies to both models.
In \cref{sec:ALM}, we briefly recall this method.
\cref{sec:POPs} is dedicated to the presentation of the results obtained
in our numerical experiment.
Therein, we consider cardinality-constrained portfolio optimization problems
from a benchmark collection and compare the performance of the augmented Lagrangian
method on \eqref{eq:CCP} and \eqref{eq:CCPref} based on function values,
iteration numbers, and the evolution of the penalty parameter.

	The particular setup of this numerical experiment has been chosen for several reasons.
	First, to the best of our knowledge, there does not exist a reasonable collection of test problems
	with vanishing constraints. 
	The instances used in \cite{HoheiselPablosPooladianSchwartzSteverango2020}
	are either toy problems or their structure is quite involved and requires detailed explanations.
	Second, a thorough numerical comparison of the bilevel optimization problem \eqref{eq:BOP} 
	and its reformulation \eqref{eq:BOPref} can be found in the recent paper \cite{ZemkohoZhou2021}
	and attests \eqref{eq:BOP} a much better behavior.
	Hence, we decided to merely investigate cardinality-constrained problems here.
	In order to ensure that the main difficulty of the test instances comes from the cardinality constraint,
	the setting of portfolio optimization has been chosen.
	Indeed, the resulting test problems are convex, quadratic optimization problems
	whenever the cardinality constraint is omitted.
	Additionally, a benchmark collection of portfolio optimization problems is available.
	Finally, to compare \eqref{eq:CCP} and its reformulation \eqref{eq:CCPref} in a fair way,
	we wanted to pick a solver which, with minor modifications, is in position to tackle both problems.
	Therefore, the augmented Lagrangian method from \cite{JiaKanzowMehlitzWachsmuth2023} has been chosen
	as it is in position to solve problems with any kind of disjunctive constraints,
	see \cite[Sections~5.1, 5.2]{JiaKanzowMehlitzWachsmuth2023},
	by employing a projected gradient method to solve the appearing subproblems.

\subsection{An augmented Lagrangian method for optimization problems with structured geometric constraints}\label{sec:ALM}

For our numerical experiment, 
we will make use of an augmented Lagrangian framework
which addresses the abstract model problem
\begin{equation}\label{eq:structured_problem}
	\min\limits_w\{\ff(w)\,|\,\gg(w)\leq 0,\,\hh(w)=0,\, w\in D\}
\end{equation}
where $\ff\colon\R^d\to\R$, $\gg\colon\R^d\to\R^{t_i}$, as well as 
$\hh\colon\R^d\to\R^{t_e}$ are continuously differentiable functions
and $D\subset\R^d$ is a nonempty, closed set which is not necessarily convex 
but allows for the fast computation of (not necessarily unique) projections.
For some penalty parameter $\rho>0$,
let us introduce the augmented Lagragian function 
$\mathcal L_\rho\colon\R^d\times\R^{t_i}\times\R^{t_e}\to\R$
of \eqref{eq:structured_problem} by means of
\begin{align*}
	&\forall w\in\R^d\,\forall \mu\in\R^{t_i}\,\forall \nu\in\R^{t_e}\colon
	\\
	&\qquad
	\mathcal L_\rho(w,\mu,\nu)
	:=
	\ff(w)+\frac{\rho}{2}\left(\norm{\max(\gg(w)+\mu/\rho,0)}^2
		+\nnorm{\hh(w)+\nu/\rho}^2\right).
\end{align*}
Note that $\mathcal L_\rho$ is a continuously differentiable function.
Furthermore, we make use of
\begin{align*}
	\forall w\in\R^d\,\forall \mu\in\R^{t_i}\colon\quad
	V_\rho(w,\mu)
	:=
	\norm{\max(\gg(w),-\mu/\rho)}+\norm{\hh(w)}
\end{align*}
in order to control the termination of the method as well as the penalty parameter.
The overall multiplier-penalty method from \cite[Section~4]{JiaKanzowMehlitzWachsmuth2023} for the numerical treatment of
\eqref{eq:structured_problem} is stated in \cref{alg:ALM}.

\begin{algorithm}
 \begin{algorithmic}[1]
        \Require $w_0 \in D $, $ \rho_0 > 0$, $\beta > 1$, $\tau \in (0,1)$, $a_{\textup{max}}\in\R^{t_i}_+$, $b_{\textup{min}}\in\R^{t_e}_-$, $b_{\textup{max}}\in\R^{t_e}_+$,
		$\{\varepsilon_k\}_{k\in\N}\subset\R_+$, $\varepsilon_{\textup{tol}}>0$
		\State{Set $\mu_0:=0$, $\nu_0:=0$, $V_0:=\infty$, and $k := 0$.}
		\While{$V_k\geq\varepsilon_{\textup{tol}}$}
		\State{Set $a_k:=\max(0,\min(\mu_k,a_{\textup{max}}))$ 
				and $b_k:=\max(b_{\textup{min}},\min(\nu_k,b_{\textup{max}}))$.}
		\parState{Compute a point $w_{k+1}\in\R^d$ such that	
    	\begin{equation*}
       		\dist(
       			-\nabla_w \mathcal{L}_{\rho_k} (w_{k+1}, a_k,b_k),
       			N_D(w_{k+1}) 
       			)
       		\leq 
       		\varepsilon_{k+1}.
    	\end{equation*}\label{item:solve_subproblem_ALM}}
		\State{Set 
			$\mu_{k+1} := \rho_k \max(\gg(w_{k+1})+a_k/\rho_k,0)$
			and
			$\nu_{k+1} := \rho_k(\hh(w_{k+1})+b_k/\rho_k)$.}
		\State{Set $V_k:=V_{\rho_k}(w_{k+1},a_k)$.}
		\If{$ k \geq 1 $ and $V_k>\tau\,V_{k-1}$}
     	\State{Set $ \rho_{k + 1} := \beta \rho_k $.} 
     	\Else
     	\State{Set $ \rho_{k + 1} := \rho_k$.}
      	\EndIf
      	\State{Set $ k := k + 1 $.}
		\EndWhile
		\State\Return $(w_k,\mu_k,\nu_k)$
	\end{algorithmic}
	\caption{Safeguarded augmented Lagrangian method for \eqref{eq:structured_problem}.}
	\label{alg:ALM}
\end{algorithm}

Apart from the appearance of the limiting normal cone to $D$ 
in \cref{item:solve_subproblem_ALM}, \cref{alg:ALM} is a rather classical
(safeguarded) augmented Lagrangian method from nonlinear programming, 
see e.g.\ \cite{BirginMartinez2014}.
Note that \cref{item:solve_subproblem_ALM} provides a requirement 
regarding the quality of the (approximate) solution
$w_{k+1}$ of the subproblem
\begin{equation}\label{eq:ALM_subproblem}
	\min\limits_w\{\mathcal L_{\rho_k}(w,a_k,b_k)\,|\,w\in D\}.
\end{equation}
Particularly, $w_{k+1}\in D$ is indispensable as one typically defines
the limiting normal cone to a closed set to be empty for out-of-set points.
By construction of the multiplier sequences 
$\{\mu_k\}_{k\in\N}$ and $\{\nu_k\}_{k\in\N}$, we can rewrite
the condition from \cref{item:solve_subproblem_ALM} by means of
\[
	\dist(
       	-\nabla_w \mathcal{L}(w_{k+1},\mu_{k+1},\nu_{k+1}),
       	N_D(w_{k+1}) 
       	)
       \leq 
    \varepsilon_{k+1},
\]
and $\mu_{k+1}\in\R^{t_i}_+$ holds by construction.
Particularly, if $\{w_{k+1}\}_{k\in\N}$ has an accumulation point $\bar w\in\R^d$, 
$\varepsilon_k\to 0$, and the multiplier sequences
are bounded, we may take the limit (along a suitable subsequence) and find
\[
	-\nabla_x\mathcal L(\bar w,\mu,\nu)\in N_D(\bar w)	
\]
for some multipliers $\mu\in\R^{t_i}_+$ and $\nu\in\R^{t_e}$.
Whenever the error measures $\{V_k\}_{k\in\N}$ tend to zero as well,
we find that $\bar w$ is feasible to \eqref{eq:structured_problem}
and $\mu$ satisfies the complementarity slackness condition w.r.t.\ the
inequality constraints.
Hence, $\bar w$ satisfies M-stationary-type conditions with multipliers $(\mu,\nu)$. 
It has been shown in \cite[Theorem~4.3, Corollary~4.4]{JiaKanzowMehlitzWachsmuth2023} 
that \cref{alg:ALM} indeed computes M-stationary points
of \eqref{eq:structured_problem} under mild assumptions.

It remains to clarify how the subproblem \eqref{eq:ALM_subproblem} can be treated 
in order to compute points fitting the requirement of \cref{item:solve_subproblem_ALM}.
For that purpose, we make use of the projected gradient method 
from \cite[Section~3]{JiaKanzowMehlitzWachsmuth2023} 
which is equipped with a nonmonotone linesearch procedure.
Recall that projections onto $D$ are easily available 
due to our standing assumptions on \eqref{eq:structured_problem}.
In the case where $\ff$ is bounded below while $D$ is bounded, 
this method does not only compute M-stationary points of the subproblem \eqref{eq:ALM_subproblem}
but the generated sequence of inner iterates comes along with validity of
the approximate M-stationarity condition demanded in \cref{item:solve_subproblem_ALM} 
for arbitrarily small $\varepsilon_{k+1}$ for sufficiently large inner iterations.

\subsection{Portfolio optimization problems with cardinality constraints}\label{sec:POPs}

For a symmetric, positive semidefinite matrix $Q\in\R^{n\times n}$, 
vectors $c,u\in\R^n_+$, and scalars $\theta>0$ as well as $\kappa\in\{1,\ldots,n-1\}$, 
we consider the cardinality-constrained portfolio optimization problem
\begin{equation}\label{eq:portfolio_cardinality_constraints}\tag{POP}
		\min\limits_z\bigl\{
			\tfrac12z^\top Qz\,\bigl|\,
			\mathtt e^\top z=1,\, c^\top z\geq\theta,\, 0\leq z\leq u,\,z\in\R^n_{\leq\kappa}
			\bigr\}.
\end{equation}
Here, the unknown vector $z$ describes the portion of $n$ individual assets 
within the portfolio which is why $\mathtt e^\top z=1$ and $z\geq 0$ appear in the constraints. 
The vector $u$ provides some upper bounds on
these individual assets while $c$ comprises expected returns for all these bonds. 
Furthermore, $\theta$ is a given lower bound for the expected return of the portfolio 
while $\kappa$ is an upper bound on the number of different bonds within it.
The matrix $Q$ describes the covariance between all the assets, so that 
\eqref{eq:portfolio_cardinality_constraints} models the problem of finding
a risk-minimal portfolio satisfying the aforementioned requirements.

We also investigate the associated complementarity-type reformulation \eqref{eq:CCPref} given by
\begin{equation}\label{eq:portfolio_cardinality_constraints_ref}\tag{POP$_\textup{ref}$}
	\min\limits_{z,\lambda}
	\left\{
		\tfrac12 z^\top Qz\,\middle|\,
		\begin{aligned}
			&\mathtt e^\top z=1,\, c^\top z\geq\theta,\, 0\leq z\leq u,\\
			&\mathtt e^\top \lambda\geq n-\kappa,\, 0\leq \lambda\leq\mathtt e,\,
			z\bullet \lambda=0
		\end{aligned}
	\right\}.
\end{equation}
Here, the difficult cardinality constraint has been replaced 
by $n$ complementarity-type constraints 
and an affine inequality on the slack variable $\lambda$,
see \cref{sec:example_CCP,sec:results_CCP} for details.

In order to solve 
\eqref{eq:portfolio_cardinality_constraints} 
and its reformulation 
\eqref{eq:portfolio_cardinality_constraints_ref} 
with the aid of \cref{alg:ALM}, 
we have to specify the role of the abstract constraint set $D$ in \eqref{eq:structured_problem}
and how to project onto it.
In our numerical experiments, the role of this abstract constraint set will be played by
\begin{equation}\label{eq:abstract_set_portfolio}
	\{z\in\R^n_{\leq\kappa}\,|\,0\leq z\leq u\}
\end{equation}
for \eqref{eq:portfolio_cardinality_constraints} and
\begin{equation}\label{eq:abstract_set_portfolio_ref}
	\{(z,\lambda)\in\R^n\times\R^n\,|\,0\leq z\leq u,\,0\leq \lambda\leq\mathtt e,\,
		z\bullet \lambda=0\}
\end{equation}
for \eqref{eq:portfolio_cardinality_constraints_ref}, respectively.
Indeed, one can easily characterize projections onto these sets, 
see \cite[Sections~5.1, 5.2]{JiaKanzowMehlitzWachsmuth2023} again.
Let us note that the set from \eqref{eq:abstract_set_portfolio} 
is the union of $\binom{n}{\kappa}$ convex polyhedral sets
while \eqref{eq:abstract_set_portfolio_ref} is the union of $2^{n}$ convex polyhedral sets.

We transferred \eqref{eq:portfolio_cardinality_constraints} and \eqref{eq:portfolio_cardinality_constraints_ref}
into the model \eqref{eq:structured_problem} by defining $D$ as stated in
\eqref{eq:abstract_set_portfolio} and \eqref{eq:abstract_set_portfolio_ref}, respectively.
Projections onto these sets are computed with the aid of the formulas
from \cite[Section~5]{JiaKanzowMehlitzWachsmuth2023}.
All the remaining constraints are augmented as indicated in \cref{sec:ALM}.
The parameters we used for the experiments are the same as in
\cite[Section~6]{JiaKanzowMehlitzWachsmuth2023}.
Furthermore, we would like to mention that \cref{alg:ALM} is automatically terminated in any of the following situations:
\begin{enumerate}
	\item\label{item:outer_iteration_counter} the (outer) iteration counter of \cref{alg:ALM} exceeds $200$,
	\item\label{item:inner_iteration_counter} the accumulated (inner) iteration counter of the projected gradient method exceeds $10^{5}$,
	\item\label{item:inverse_stepsize} the (inverse) stepsize used in the projected gradient method falls below $10^{-20}$, or
	\item\label{item:penalty_parameter} the penalty parameter exceeds $10^{18}$.
\end{enumerate}

For our numerical experiments, 
we exploited the test problem collection used in \cite{FrangioniGentile2007},
which is available from the webpage
    \url{https://commalab.di.unipi.it/datasets/MV/}.
Here, we used all 30 test instances of dimension $n:=200$ and the three different
cardinality bounds $\kappa\in\{5,10,20\}$ for each problem. 
In \cref{tab:convex_branches}, we illustrate the numbers 
of convex branches within the sets from
\eqref{eq:abstract_set_portfolio} and \eqref{eq:abstract_set_portfolio_ref}.
\begin{table}
	\centering
	\begin{tabular}{c c c c}
	\toprule
	$\binom{200}{5}$	&	$\binom{200}{10}$	&	$\binom{200}{20}$	&	$2^{200}$
	\\
	\midrule
	$2.5\cdot 10^{9}$	&	$2.2\cdot 10^{17}$	&	$1.6\cdot 10^{27}$	&	$1.6\cdot 10^{60}$
	\\
	\bottomrule
	\end{tabular}
	\caption{Number of convex branches in the sets from \eqref{eq:abstract_set_portfolio}
		and \eqref{eq:abstract_set_portfolio_ref}.}
	\label{tab:convex_branches}
\end{table}

\begin{figure}[htp]
\centering
	\pgfplotstableread[col sep=tab]{portfolio_alm_val.txt}\almtable
	\pgfplotstableread[col sep=tab]{portfolio_alm_ref_val.txt}\almreftable
	\begin{tikzpicture}
  	\begin{axis}[
      width=1.0\textwidth,
      height=6cm,
      ybar=0cm, 
      bar width=.1cm,
      enlarge x limits={abs=.2cm},
      ylabel={objective function value},
      xtick=data,     
      xticklabels from table={\almtable}{Name of problem},
      x tick label style={ rotate=45,scale=.7,anchor=east }
    ]
    \addplot [fill=red!90!black] 
    	table [x expr=\coordindex, y={kappa 20}] {\almtable};
    \addplot [fill=yellow!90!black] 
    	table [x expr=\coordindex, y={kappa 20}] {\almreftable};
  	\end{axis}
	\end{tikzpicture}
	\\
	 \begin{tikzpicture}
  	\begin{axis}[
      width=1.0\textwidth,
      height=6cm,
      ybar=0cm, 
      bar width=.1cm,
      enlarge x limits={abs=.2cm},
      ylabel={objective function value},
      xtick=data,     
      xticklabels from table={\almtable}{Name of problem},
      x tick label style={ rotate=45,scale=.7,anchor=east }
    ]
    \addplot [fill=red!90!black] 
    	table [x expr=\coordindex, y={kappa 10}] {\almtable};
    \addplot [fill=yellow!90!black] 
    	table [x expr=\coordindex, y={kappa 10}] {\almreftable};
  	\end{axis}
	\end{tikzpicture}
	\\
	\begin{tikzpicture}
  	\begin{axis}[
      width=1.0\textwidth,
      height=6cm,
      ybar=0cm, 
      bar width=.1cm,
      enlarge x limits={abs=.2cm},
      ylabel={objective function value},
      xtick=data,     
      xticklabels from table={\almtable}{Name of problem},
      x tick label style={ rotate=45,scale=.7,anchor=east }
    ]
    \addplot [fill=red!90!black] 
    	table [x expr=\coordindex, y={kappa 5}] {\almtable};
    \addplot [fill=yellow!90!black] 
    	table [x expr=\coordindex, y={kappa 5}] {\almreftable};
  	\end{axis}
	\end{tikzpicture}
	\caption{Final function values obtained by \cref{alg:ALM} when
		applied to the original problem \eqref{eq:portfolio_cardinality_constraints} (red)
		and its reformulation \eqref{eq:portfolio_cardinality_constraints_ref} (yellow) 
		with cardinality levels
		$\kappa = 20$, $\kappa=10$, and $ \kappa = 5 $ (top to bottom), respectively.}
	\label{fig:SparseKappa_val}
\end{figure}
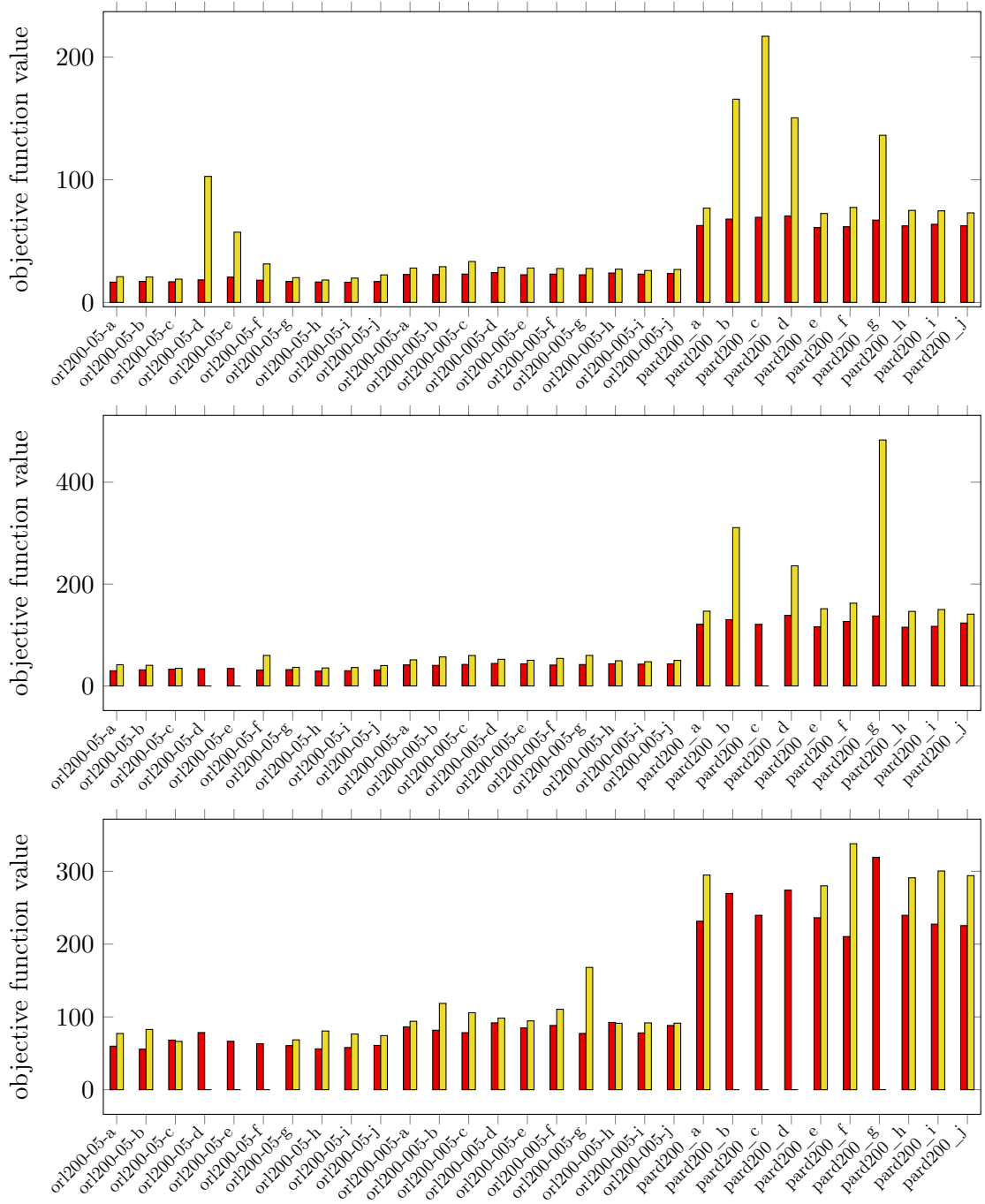

In order to solve \eqref{eq:portfolio_cardinality_constraints}, 
we always used $w^0:=z^0:=0$ as a starting point.
For \eqref{eq:portfolio_cardinality_constraints_ref}, we used $z^0:=0$ and
$\lambda^0\in\R^n$ given by $\lambda^0_i:=1$ for $i=1,\ldots,n-\kappa$ 
and $\lambda^0_i:=0$ for $i=n-\kappa+1,\ldots,n$
in order to construct the starting point $w^0:=(z^0,\lambda^0)$.
Note that these points belong to \eqref{eq:abstract_set_portfolio} 
and \eqref{eq:abstract_set_portfolio_ref}, respectively.
The results of our experiments are depicted in \cref{fig:SparseKappa_val}.
Therein, the objective function values of the final iterates of \cref{alg:ALM} are plotted.
Note that whenever \cref{alg:ALM} is terminated due to 
one of the four reasons~\ref{item:outer_iteration_counter}-\ref{item:penalty_parameter} 
mentioned above, a bar of zero height is plotted in \cref{fig:SparseKappa_val}. 

Only for 2 of the 90 test problems, the better objective value is generated via the
reformulated problem \eqref{eq:portfolio_cardinality_constraints_ref}. In all
other settings, the direct use of the original model \eqref{eq:portfolio_cardinality_constraints}
turned out to be beneficial. Furthermore, we would like to mention that,
when applied to \eqref{eq:portfolio_cardinality_constraints_ref}, \cref{alg:ALM} has been
aborted due to one of the reasons~\ref{item:outer_iteration_counter}-\ref{item:penalty_parameter}
for 10 of the test problems having cardinality level $\kappa\in\{5,10\}$. Here, termination 
was caused by~\ref{item:inner_iteration_counter} (for 5 test problems) 
and~\ref{item:inverse_stepsize} (for 5 test problems).

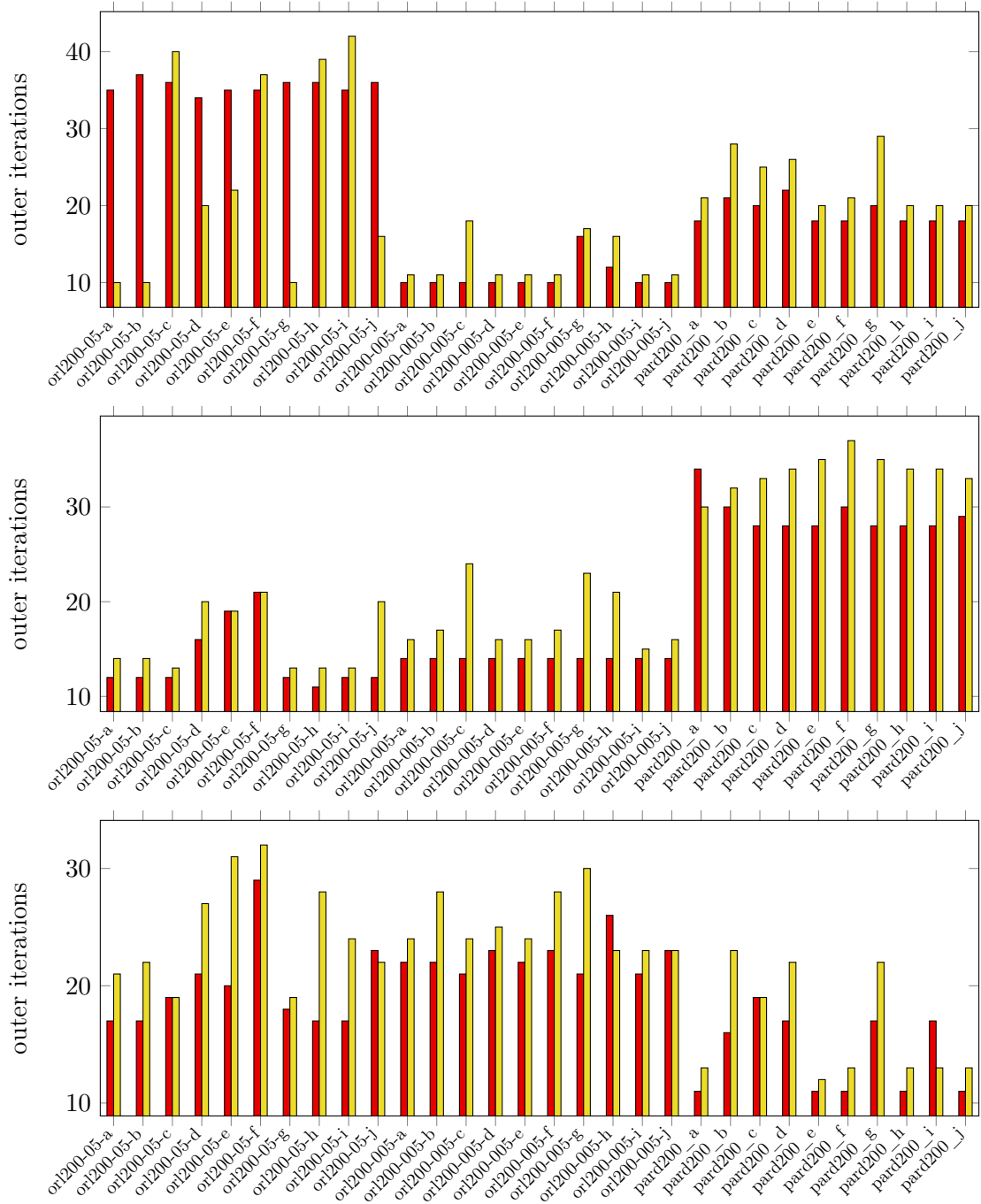
\begin{figure}[htp]
\centering
	\pgfplotstableread[col sep=tab]{portfolio_alm_outer_iter.txt}\almtable
	\pgfplotstableread[col sep=tab]{portfolio_alm_ref_outer_iter.txt}\almreftable
	\begin{tikzpicture}
  	\begin{axis}[
      width=1.0\textwidth,
      height=6cm,
      ybar=0cm, 
      bar width=.1cm,
      enlarge x limits={abs=.2cm},
      ylabel={outer iterations},
      xtick=data,     
      xticklabels from table={\almtable}{Name of problem},
      x tick label style={ rotate=45,scale=.7,anchor=east }
    ]
    \addplot [fill=red!90!black] 
    	table [x expr=\coordindex, y={kappa 20}] {\almtable};
    \addplot [fill=yellow!90!black] 
    	table [x expr=\coordindex, y={kappa 20}] {\almreftable};
  	\end{axis}
	\end{tikzpicture}
	\\
	 \begin{tikzpicture}
  	\begin{axis}[
      width=1.0\textwidth,
      height=6cm,
      ybar=0cm, 
      bar width=.1cm,
      enlarge x limits={abs=.2cm},
      ylabel={outer iterations},
      xtick=data,     
      xticklabels from table={\almtable}{Name of problem},
      x tick label style={ rotate=45,scale=.7,anchor=east }
    ]
    \addplot [fill=red!90!black] 
    	table [x expr=\coordindex, y={kappa 10}] {\almtable};
    \addplot [fill=yellow!90!black] 
    	table [x expr=\coordindex, y={kappa 10}] {\almreftable};
  	\end{axis}
	\end{tikzpicture}
	\\
	\begin{tikzpicture}
  	\begin{axis}[
      width=1.0\textwidth,
      height=6cm,
      ybar=0cm, 
      bar width=.1cm,
      enlarge x limits={abs=.2cm},
      ylabel={outer iterations},
      xtick=data,     
      xticklabels from table={\almtable}{Name of problem},
      x tick label style={ rotate=45,scale=.7,anchor=east }
    ]
    \addplot [fill=red!90!black] 
    	table [x expr=\coordindex, y={kappa 5}] {\almtable};
    \addplot [fill=yellow!90!black] 
    	table [x expr=\coordindex, y={kappa 5}] {\almreftable};
  	\end{axis}
	\end{tikzpicture}
	\caption{Number of (outer) iterations done by \cref{alg:ALM} when
		applied to the original problem \eqref{eq:portfolio_cardinality_constraints} (red)
		and its reformulation \eqref{eq:portfolio_cardinality_constraints_ref} (yellow) with cardinality 
		levels $\kappa = 20$, $\kappa=10$, and $ \kappa = 5 $ (top to bottom), respectively.}
	\label{fig:SparseKappa_outer_iter}
\end{figure}

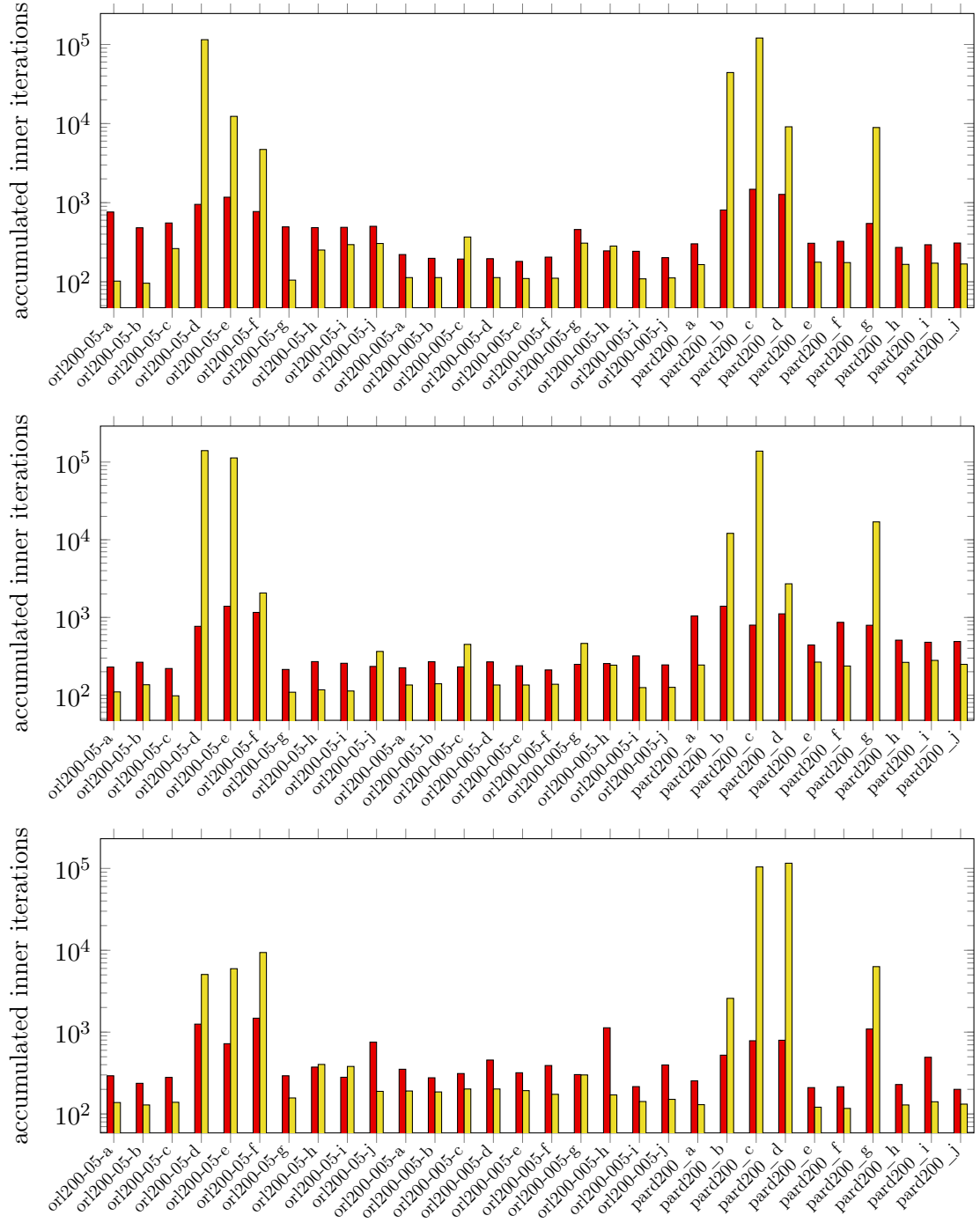
\begin{figure}[htp]
\centering
	\pgfplotstableread[col sep=tab]{portfolio_alm_inner_iter.txt}\almtable
	\pgfplotstableread[col sep=tab]{portfolio_alm_ref_inner_iter.txt}\almreftable
	\begin{tikzpicture}
  	\begin{axis}[
      width=1.0\textwidth,
      height=6cm,
      ybar=0cm, 
      bar width=.1cm,
      enlarge x limits={abs=.2cm},
      ylabel={accumulated inner iterations},
      xtick=data,     
      xticklabels from table={\almtable}{Name of problem},
      x tick label style={ rotate=45,scale=.7,anchor=east },
      ymode=log
    ]
    \addplot [fill=red!90!black] 
    	table [x expr=\coordindex, y={kappa 20}] {\almtable};
    \addplot [fill=yellow!90!black] 
    	table [x expr=\coordindex, y={kappa 20}] {\almreftable};
  	\end{axis}
	\end{tikzpicture}
	\\
	 \begin{tikzpicture}
  	\begin{axis}[
      width=1.0\textwidth,
      height=6cm,
      ybar=0cm, 
      bar width=.1cm,
      enlarge x limits={abs=.2cm},
      ylabel={accumulated inner iterations},
      xtick=data,     
      xticklabels from table={\almtable}{Name of problem},
      x tick label style={ rotate=45,scale=.7,anchor=east },
      ymode=log
    ]
    \addplot [fill=red!90!black] 
    	table [x expr=\coordindex, y={kappa 10}] {\almtable};
    \addplot [fill=yellow!90!black] 
    	table [x expr=\coordindex, y={kappa 10}] {\almreftable};
  	\end{axis}
	\end{tikzpicture}
	\\
	\begin{tikzpicture}
  	\begin{axis}[
      width=1.0\textwidth,
      height=6cm,
      ybar=0cm, 
      bar width=.1cm,
      enlarge x limits={abs=.2cm},
      ylabel={accumulated inner iterations},
      xtick=data,     
      xticklabels from table={\almtable}{Name of problem},
      x tick label style={ rotate=45,scale=.7,anchor=east },
      ymode=log
    ]
    \addplot [fill=red!90!black] 
    	table [x expr=\coordindex, y={kappa 5}] {\almtable};
    \addplot [fill=yellow!90!black] 
    	table [x expr=\coordindex, y={kappa 5}] {\almreftable};
  	\end{axis}
	\end{tikzpicture}
	\caption{Number of accumulated inner iterations of the subproblem solver
		when \cref{alg:ALM} is
		applied to the original problem \eqref{eq:portfolio_cardinality_constraints} (red)
		and its reformulation \eqref{eq:portfolio_cardinality_constraints_ref} (yellow) with cardinality 
		levels $\kappa = 20$, $\kappa=10$, and $ \kappa = 5 $ (top to bottom), respectively.}
	\label{fig:SparseKappa_inner_iter}
\end{figure}

\begin{figure}[htp]
\centering
	\pgfplotstableread[col sep=tab]{portfolio_alm_pp.txt}\almtable
	\pgfplotstableread[col sep=tab]{portfolio_alm_ref_pp.txt}\almreftable
	\begin{tikzpicture}
  	\begin{axis}[
      width=1.0\textwidth,
      height=6cm,
      ybar=0cm, 
      bar width=.1cm,
      enlarge x limits={abs=.2cm},
      ylabel={final penalty parameter},
      xtick=data,     
      xticklabels from table={\almtable}{Name of problem},
      x tick label style={ rotate=45,scale=.7,anchor=east },
      ymode=log
    ]
    \addplot [fill=red!90!black] 
    	table [x expr=\coordindex, y={kappa 20}] {\almtable};
    \addplot [fill=yellow!90!black] 
    	table [x expr=\coordindex, y={kappa 20}] {\almreftable};
  	\end{axis}
	\end{tikzpicture}
	\\
	 \begin{tikzpicture}
  	\begin{axis}[
      width=1.0\textwidth,
      height=6cm,
      ybar=0cm, 
      bar width=.1cm,
      enlarge x limits={abs=.2cm},
      ylabel={final penalty parameter},
      xtick=data,     
      xticklabels from table={\almtable}{Name of problem},
      x tick label style={ rotate=45,scale=.7,anchor=east },
      ymode=log
    ]
    \addplot [fill=red!90!black] 
    	table [x expr=\coordindex, y={kappa 10}] {\almtable};
    \addplot [fill=yellow!90!black] 
    	table [x expr=\coordindex, y={kappa 10}] {\almreftable};
  	\end{axis}
	\end{tikzpicture}
	\\
	\begin{tikzpicture}
  	\begin{axis}[
      width=1.0\textwidth,
      height=6cm,
      ybar=0cm, 
      bar width=.1cm,
      enlarge x limits={abs=.2cm},
      ylabel={final penalty parameter},
      xtick=data,     
      xticklabels from table={\almtable}{Name of problem},
      x tick label style={ rotate=45,scale=.7,anchor=east },
      ymode=log
    ]
    \addplot [fill=red!90!black] 
    	table [x expr=\coordindex, y={kappa 5}] {\almtable};
    \addplot [fill=yellow!90!black] 
    	table [x expr=\coordindex, y={kappa 5}] {\almreftable};
  	\end{axis}
	\end{tikzpicture}
	\caption{Final value of the penalty parameter in \cref{alg:ALM} when
		applied to the original problem \eqref{eq:portfolio_cardinality_constraints} (red)
		and its reformulation \eqref{eq:portfolio_cardinality_constraints_ref} (yellow) with cardinality 
		levels $\kappa = 20$, $\kappa=10$, and $ \kappa = 5 $ (top to bottom), respectively.}
	\label{fig:SparseKappa_pp}
\end{figure}
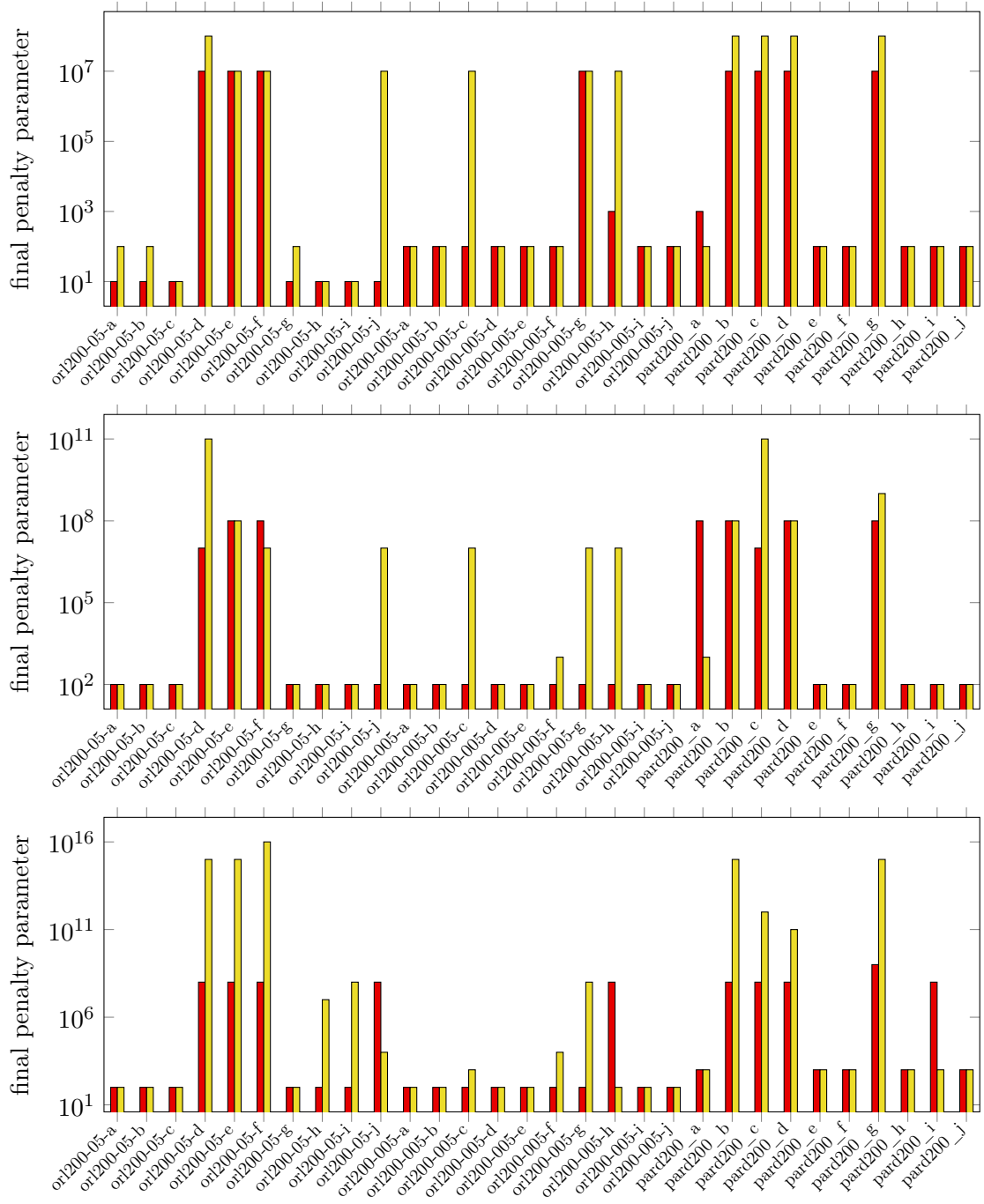

In \cref{fig:SparseKappa_outer_iter,fig:SparseKappa_inner_iter,fig:SparseKappa_pp}, 
we use a similar
type of visualization in order to document the number of (outer) iterations in \cref{alg:ALM},
the number of accumulated inner iterations done by the projected gradient method in order to
solve the appearing subproblems, and the final value of the penalty parameter.
For most of the test problems, \cref{alg:ALM} needed less outer but more inner iterations when
applied to the original model \eqref{eq:portfolio_cardinality_constraints}, see
\cref{fig:SparseKappa_outer_iter,fig:SparseKappa_inner_iter}. This might be caused
by the fact that the disjunctive structure of the set from \eqref{eq:abstract_set_portfolio} is
more complicated than the one of \eqref{eq:abstract_set_portfolio_ref}
although the number of its convex branches is significantly smaller, 
see \cref{tab:convex_branches} again.
Furthermore, for most of the test problems, 
model \eqref{eq:portfolio_cardinality_constraints} and \eqref{eq:portfolio_cardinality_constraints_ref}
came along with the same value of the final penalty parameter.
In the case where these values were different,
the smaller one has been realized via \eqref{eq:portfolio_cardinality_constraints} more frequently.

These results clearly indicate that solving \eqref{eq:portfolio_cardinality_constraints} via
\cref{alg:ALM} directly is superior to the treatment of the reformulated problem
\eqref{eq:portfolio_cardinality_constraints_ref} with the same method. In this regard,
our results clearly advise against making the implicit variable $\lambda$ explicit.

\section{Concluding remarks}\label{sec:conclusions}

In this paper,
we have shown in terms of the model problem \eqref{eq:implicit_problem} and its reformulation \eqref{eq:explicit_problem}
that treating implicit variables as explicit ones is likely to be disadvantageous
from both a theoretical and computational point of view.
Indeed, it turned out that \eqref{eq:explicit_problem} is likely to possess more
stationary points (of different types) and, particularly, local minimizers than \eqref{eq:implicit_problem}.
Furthermore, we have shown that comparable constraint qualifications which address \eqref{eq:explicit_problem}
are often more restrictive than their associated counterparts which suit \eqref{eq:implicit_problem}.
In this regard, we recover and even deepen the insights obtained in \cite{BenkoMehlitz2021}
on the base of a much simpler model problem which is easy to access.
Furthermore, it has been illustrated by means of a numerical experiment that exploiting 
\eqref{eq:explicit_problem} for numerical purposes may turn out to be worse than simply
tackling \eqref{eq:implicit_problem}. These findings are in line with the computational studies
from \cite{Mehlitz2020c} and \cite{ZemkohoZhou2021} which address or-constrained and bilevel optimization, respectively.

\subsection*{Acknowledgments}

The author acknowledges a fruitful discussion with Mat\'{u}\v{s} Benko
on an earlier version of this paper.
Additionally, the remarks of a reviewer helped to improve the presentation 
of the obtained results.

\end{document}